\documentclass[a4wide]{article}
\usepackage[utf8]{inputenc}
\usepackage[a4paper, total={6.5in, 9.25in}]{geometry}

\usepackage{hyperref}
\usepackage[auth-lg]{authblk}
\usepackage{amsfonts}
\usepackage{babel}[english]

\usepackage{bm}
\usepackage{mathtools}
\usepackage{eufrak}

\usepackage{enumerate}

\usepackage{algorithm,algorithmic}

\usepackage{caption}
\usepackage{subcaption}

\usepackage{xcolor}
 
\usepackage{diagbox}
\usepackage{bigints}
\usepackage{relsize}

\DeclarePairedDelimiter\floor{\lfloor}{\rfloor}
\DeclareRobustCommand{\rchi}{{\mathpalette\irchi\relax}}
\newcommand{\irchi}[2]{\raisebox{\depth}{$#1\chi$}}

\usepackage[page]{appendix} 

\usepackage{comment}

\newtheorem{remark}{Remark}

\begin{document}

\newcommand{\WidestEntry}{$\left[13,8,555 \right]$}
\newcommand{\SetToWidest}[1]{\makebox[\widthof{\WidestEntry}]{#1}}

\title{Efficient approximation of cardiac mechanics through reduced order modeling with deep learning-based operator approximation}


\author[1]{Ludovica Cicci}
\author[1]{Stefania Fresca}
\author[1]{Andrea Manzoni}
\author[1,2]{Alfio Quarteroni}
\setlength{\affilsep}{1em}
\renewcommand\Authsep{, }
\affil[1]{MOX-Dipartimento di Matematica, Politecnico di Milano, P.zza Leonardo da Vinci 32, 20133 Milano, Italy\footnote{\texttt{{ludovica.cicci,stefania.fresca,andrea1.manzoni,alfio.quarteroni}@polimi.it}}}
\affil[2]{Mathematics Institute, \'Ecole Polytechnique F\'ed\'erale de Lausanne, Switzerland (Professor Emeritus)}
\date{}
\renewcommand\Affilfont{\small}

\maketitle

\section*{Abstract}

Reducing the computational time required by high-fidelity, full order models (FOMs) for the solution of problems in cardiac mechanics is crucial to allow the translation of patient-specific simulations into clinical practice. While FOMs, such as those based on the finite element method, provide valuable information of the cardiac mechanical function, up to hundreds of thousands degrees of freedom may be needed to obtain accurate numerical results. As a matter of fact, simulating even just a few heartbeats can require hours to days of CPU time even on powerful supercomputers. In addition, cardiac models depend on a set of input parameters that we could let vary in order to explore multiple virtual scenarios. To compute reliable solutions at a greatly reduced computational cost, we rely on a reduced basis method empowered with a new deep-learning based operator approximation, which we refer to as \textit{Deep-HyROMnet} technique. Our strategy combines a projection-based POD-Galerkin method with deep neural networks for the approximation of (reduced) nonlinear operators, overcoming the typical computational bottleneck associated with standard hyper-reduction techniques. This method is shown to provide reliable approximations to cardiac mechanics problems outperforming classical projection-based ROMs in terms of computational speed-up of orders of magnitude, and enhancing forward uncertainty quantification analysis otherwise unaffordable.


\vspace{0.75cm}

\section{Introduction}

Cardiac mechanics simulations aim at reproducing the response of the cardiac muscle under specified loading conditions and consist of large-scale differential systems governed by the equations of elasto\-dynamics, complemented with suitable constitutive laws to correctly capture the mechanical behavior of the myocardium. Modeling the cardiac dynamics is therefore a challenging task, as the myocardium is a strongly non-isotropic, incompressible material, characterized by an exponential nonlinear strain energy function \cite{guccione1995finite, holzapfel2009constitutive} and a fiber-sheet structure \cite{streeter1969fiber}. Another crucial aspect is the inclusion in the mathematical model of the active forces that drive the contraction mechanism of the muscle, which is able to contract after being electrically activated; these active properties are time-dependent and anisotropic. 

In the last decades, there has been substantial advances in the development of multi-physics, multi-scale mathematical models of cardiac functions \cite{kerckhoffs2007coupling, goktepe2010electromechanics, trayanova2011whole, gurev2015high, niederer2019computational, regazzoni2020cardiac}. The availability of realistic, patient-specific simulations, both for normal and diseased hearts, allows for a quantitative understanding of cardiac physiology and raises the prospect of their use in a number of applications, e.g., for improving diagnosis, providing real-time decision support, predicting prognosis and, ultimately, supporting clinical decisions  \cite{winslow2012computational, niederer2016using, chabiniok2016multiphysics, mangion2018advances}. However, the translation of cardiac simulations into the clinical practice is often hampered by the huge computational costs involved in the solution to the underlying problem by means of suitable numerical procedures, such as the finite element method (FEM) \cite{guccione1991finite, nash2000computational, lafortune2012coupled, quarteroni2017cardiovascular, mcculloch2020large}.
 
While finite element models of the heart provide valuable information, they may require up to hundreds of thousands degrees of freedom (dofs) to obtain accurate numerical results, so that simulating even just a few heartbeats can require hours to days of CPU time even on supercomputers. Additionally, cardiac models depend on a large set of patient-specific parameters characterizing, e.g., material prop\-er\-ties, boundary/initial conditions, geometrical features, or local fiber orientation, which are affected by un\-cer\-tain\-ty and should be properly calibrated through optimization routines. Being able to perform efficient numerical simulations in this context is indeed essential to explore multiple virtual scenarios, to quantify cardiac outputs and related uncertainties, as well as to evaluate the impact of pathological conditions. All these tasks require repeated model evaluations over different input parameter values, thus making high-fidelity, full order models (FOMs) computationally unaffordable. 

Alternative numerical methods have been developed in the past decades aiming to compute reliable solutions to parameter-dependent problems at a greatly reduced computational cost, such as data-driven surrogate models and projection-based reduced order models (ROMs). The former aim to learn, in a non-intrusive way, the hidden relation between input parameters and corresponding output quantities of interest (possibly including problem's solution) from usually large data sets of input-output pairs. For instance, in \cite{rodriguez2019uncertainty, campos2020uncertainty} surrogate models were generated via the polynomial chaos expansion approach to accelerate uncertainty quantification (UQ) studies and sensitivity analysis of left ventricular mechanics. Many machine learning-based models have been proposed as real-time cardiac mechanics simulators \cite{dabiri2019prediction, maso2020deep, dalton2021graph}, while statistical emulators, such as Gaussian processes, have been used to speed-up parameter inference \cite{borowska2020bayesian, noe2019gaussian} or to reduce the complexity of parametric searches for high-fidelity models \cite{di2018gaussian}. Despite being well suited for the rapid and repeated evaluation of the input-output map, this models may lack of accuracy when dealing with a single, patient-specific forward simulation of the cardiac activity. Moreover, the numerical test cases presented are mostly restricted to the diastolic filling phase, when only passive material properties are taken into account, and to the solution of quasi-static mechanics problems. 

On the contrary, projection-based ROMs, such as the Galerkin-reduced basis (RB) method, replace the high-fidelity problem with a reduced problem featuring lower computational complexity, still retaining the essential features of the FOM. These methods are usually characterized by a splitting of the reduction procedure into an expensive offline phase, during which multiple parametric instances of the FOM are computed to generate a basis for the reduced subspace, and an efficient online phase. A reduction strategy for the quasi-static mechanics problem is proposed in \cite{bonomi2017reduced}, where proper orthogonal decomposition (POD) for basis construction is combined with suitable hyper-reduction techniques to efficiently handle nonlinear terms, whereas in \cite{cicci2021cardiacDEIM} POD-Galerkin ROMs exploiting the discrete empirical interpolation method (DEIM) have been exploited for the efficient and accurate solution to the time-dependent cardiac problem, on both idealized and patient-specific left ventricle geometries, albeit using a relative low number of degrees of freedom. In \cite{hirschvogel2018computational, pfaller2020using} POD has been applied to reduce the structural dimension of a monolithic 3D-0D coupled structure-circulation model in a four-chamber, patient-specific geometry. Nonetheless, despite their application in a wide range of scenarios, relatively contained speed-ups are achieved by projection-based ROMs in cardiac mechanics due to the highly nonlinear nature of the problem. Indeed, if the construction of a reduced subspace to approximate the problem solution does not pose serious issues, resulting in extremely low dimensional spaces even for complex material laws, the bottleneck in all these cases is represented by the assembling of reduced operators, and the projection of the approximated operators through DEIM.  \\

Motivated by this observation, in this work we address the efficient solution to parameterized cardiac mechanics problems by means of our newly developed \textit{Deep-HyROMnet} method \cite{cicci2021DeepHyROMnet}. The key idea of this projection- and deep-learning-based method is to exploit within the Galerkin-RB approach suitable deep neural network (DNN) architectures -- as the ones introduced in \cite{fresca2021comprehensive,fresca2021pod} -- to approximate reduced nonlinear operators efficiently. Unlike data-driven strategies, for which the predicted output is not guaranteed to satisfy the underlying PDE, Deep-HyROMnet is a fully physics-based ROM, as it computes the problem solution by solving a reduced nonlinear system built by enforcing the problem's equations onto a (linear) reduced-order subspace. In this work, we show how Deep-HyROMnet outperforms classical POD-Galerkin-DEIM ROMs in terms of computational speed-up for the solution to a 3D-0D coupled structure-circulation model for the left ventricle, both in physiological and pathological scenarios. By providing accurate and computationally efficient simulations of the left ventricle dynamics, the reduction strategy is successfully used to address the solution to many-query tasks, specifically  forward UQ. \\

The reminder of the paper is structured as follows. After a brief introduction of the basic concepts of continuum mechanics, in Sec.~\ref{sec:3D-0D} we provide the 3D elastodynamics model for the cardiac tissue and the 0D hemodynamics model for the blood circulation. Further, the high-fidelity formulation of the monolithically coupled structure-circula\-tion model for the description of the mechanical activity of the left ventricle during a whole heartbeat is presented. In Sec.~\ref{sec:LV-Deep-HyROMnet} we show how Deep-HyROMnet can be employed in this context, whilst the numerical performances of the resulting hyper-reduced ROM are assessed in Sec.~\ref{sec:numerical_results} on two different applications, the former focusing on a physiological scenario and the latter assuming the presence of an ischemic region inside the myocardium. Preliminary results on the application of Deep-HyROMnet in the multi-query context of forward UQ are then presented in Sec.~\ref{sec:forwardUQ}, whilst conclusions are drawn in Sec.~\ref{sec:conclusion}.
Details on the Deep-HyROMnet technique, as well as the POD-Galerkin-DEIM method used as benchmark model, are reported in the Appendices~\ref{appendix:POD} and \ref{appendix:DeepHyROMnet} to make the paper self-contained.


\section{Mathematical models: 3D-0D mechanics-circulation model}\label{sec:3D-0D}

\begin{figure}[b!]
	\centering
	\includegraphics[width=0.9\textwidth]{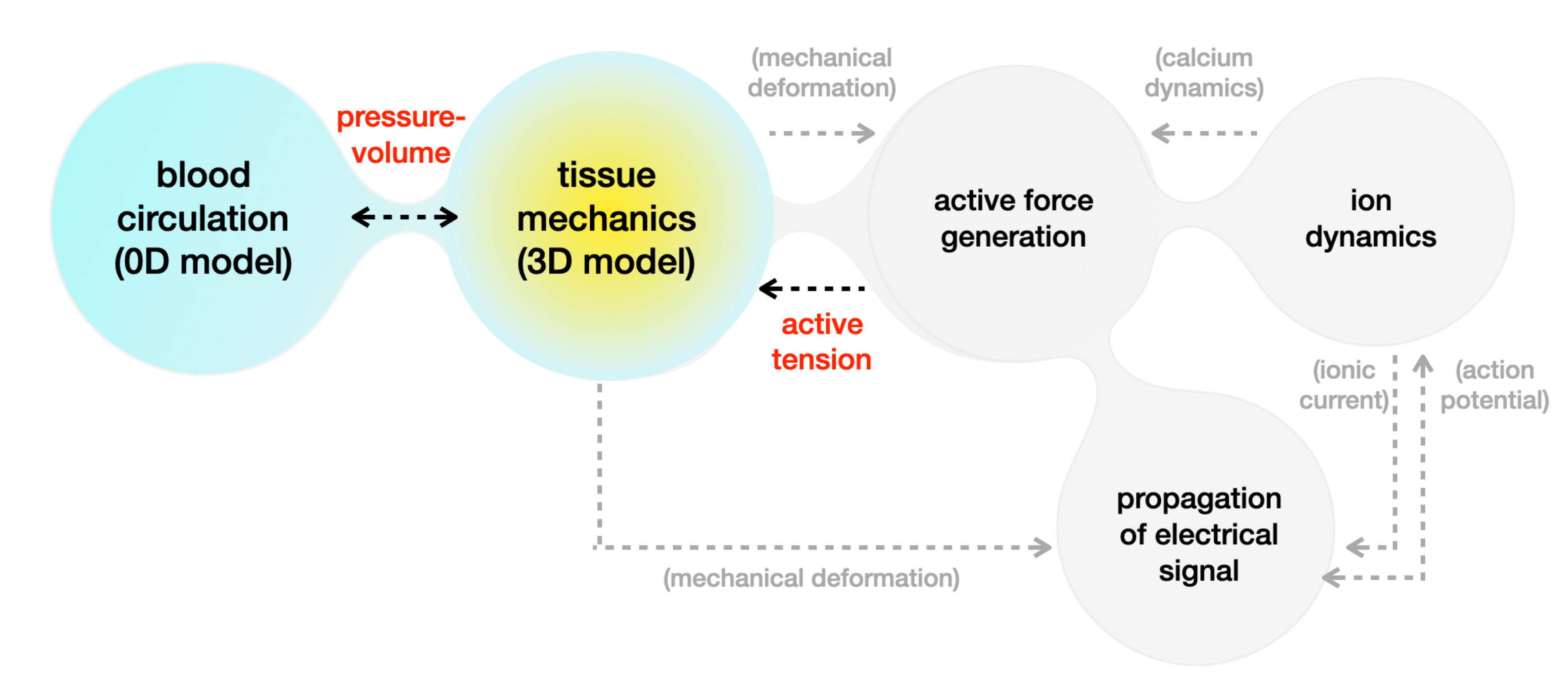}
			\caption{Cardiac core models: the two highlighted blocks on the left represent the two 3D and 0D submodels yielding the coupled mechanics-circulation model addressed in this work. The grey blocks on the right (active force generation, ion dynamics, and propagation of electrical signal) are not addressed by the proposed framework; their response is surrogated by the definition of a suitable, analytical, function to encode the behavior of the active function, that acts as a physical coefficient for the tissue mechanics (3D) model. Interactions among the submodels reported in grey are also discarded.}
	\label{fig:core_models}
\end{figure}

The solution of cardiac mechanics problems involves the interaction between several biophysical phenomena concurring to the heart function, namely electro\-physiology, biochemistry, mechanics and fluid dynamics, each described by suitable models (see Figure~\ref{fig:core_models}) written in terms of PDEs and/or ODEs \cite{quarteroni2017cardiovascular}. Ele\-ctro\-phy\-siolo\-gy corresponds to the propagation of the electrical potential and ion dynamics, and describes the electric activity of cardiac muscle cells; the activation of cardiomyocytes is the results of complex mechano-chemical interactions among contractile proteins \cite{bers2001excitation} and provides the active tension necessary to the mechanics model. In this work, we focus on the mechanical behavior only, surrogating the active force generation model through an explicit, periodic, analytical active tension function. See, e.g., \cite{gerbi2019monolithic, regazzoni2020machine, regazzoni2020cardiac} for a detailed presentation of a fully coupled cardiac electromechanics model. 

For the sake of completeness, in Sec.~\ref{sec:3D} we report the three-dimensional (3D) model used to describe the mechanical activity of the left ventricle during the cardiac cycle, taking into account both the passive response of the tissue due to the presence of blood and the active behavior of the muscular fibers. The structural model is then coupled to a lumped-parameter fluid zero-dimensional (0D) model  -- that only involves a system of ordinary differential equations, thus depending on the time variable only -- to provide the unknown pressure to the endocardial wall (see Sec.~\ref{sec:0D}). Since biological tissues commonly experience large deformations \cite{demiray1981large}, we rely on the finite elasticity theory for the description of cardiac mechanics by means of nonlinear, time-dependent partial differential equations (PDEs).


\subsection{3D elastodynamics model}\label{sec:3D}

Given a continuum body $\mathcal{B}$ embedded in a three-dimensional Euclidean space, let $\Omega_0\subset\mathbb{R}^3$ be its reference configuration at time $t=0$ and $\Omega_t\subset\mathbb{R}^3$ its current configuration at time $t>0$. The motion of the body $\mathcal{B}$ is defined as the map $\chi\colon\Omega_0\rightarrow\Omega_t$, for all $t>0$, such that $\mathbf{x} = \chi(\mathbf{X},t)$, where $\mathbf{X}$ and $\mathbf{x}$ denote the position vectors in the reference $\Omega_0$ and in the current $\Omega_t$ configurations, respectively. The effects of deformation on the solid body $\mathcal{B}$ are described by means of the displacement field
\begin{equation*}
\mathbf{u}(\mathbf{X},t;\bm\mu) = \chi(\mathbf{X},t) - \mathbf{X},
\end{equation*}
which represents the unknown of our problem and depends on a set of model parameters, such as material coefficients, boundary/initial conditions, source terms and so on, collectively denoted by $\bm\mu\in\mathcal{P}\subset\mathbb{R}^P$, where $\mathcal{P}$ is a compact set. Other important quantities in the framework of continuum mechanics are the deformation gradient $\mathbf{F}$, the right Cauchy-Green tensor $\mathbf{C}$ and the Green-Lagrange strain tensor $\mathbf{E}$, that are defined as
\begin{equation*}
\mathbf{F} = \frac{\partial \mathbf{x}}{\partial \mathbf{X}}, \quad \mathbf{C} = \mathbf{F}^T\mathbf{F} \text{~~and~~} \mathbf{E} = \frac{1}{2}(\mathbf{C} - \mathbf{I}),
\end{equation*}
respectively. The change in volume between the reference and the current configuration at time $t>0$ is given by the determinant of the deformation gradient, i.e. $J(\mathbf{X},t) = \det\mathbf{F}(\mathbf{X},t)$, known as the volume ratio. A motion for which $J=1$ is said to be \textit{isochoric} or \textit{isovolumetric}.

The displacement field $\mathbf{u}(\bm\mu)\in\Omega_0\times\mathbb{R}^+$, for $\bm\mu\in\mathcal{P}$, can be found by solving the equation of motion given by the balance of linear momentum~\cite{Holzapfel}, that is
\begin{equation*}\label{eq:motionCauchy}
\rho_0\partial_t^2\mathbf{u}(\mathbf{X},t;\bm\mu) - \nabla_0\cdot\mathbf{P}(\mathbf{F}(\mathbf{X},t;\bm\mu)) = \mathbf{b}_0(\mathbf{X},t;\bm\mu), \qquad \mathbf{X}\in\Omega_0, t>0,
\end{equation*}
when boundary and initial conditions are provided. Here, $\rho_0>0$ denotes the density of the body, $\mathbf{P}$ is the first Piola-Kirchhoff stress tensor and $\mathbf{b}_0$ represents a body force field. Suitable constitutive law, i.e. stress-strain relationships, must be specified  to describe the behavior of the given material. Furthermore, to incorporate active contraction of the tissue, we adopt an active stress approach \cite{ambrosi2012active}, which assumes an additive decomposition of the stress tensor into a passive and an active contributions as
\begin{equation*}
\begin{aligned}
\mathbf{P}(\mathbf{F}) = \mathbf{P}_p(\mathbf{F}) + \mathbf{P}_a(\mathbf{F}).
\end{aligned}
\end{equation*}
For what concerns the passive term $\mathbf{P}_p(\mathbf{F})$, we consider the myocardium as hyperelastic, for which we can assume the existence of a strain density function $\mathcal{W}\colon Lin^+\rightarrow\mathbb{R}$ such that
\begin{equation*}
\mathbf{P}_p(\mathbf{F}) = \frac{\partial\mathcal{W}(\mathbf{F})}{\partial\mathbf{F}}.
\end{equation*}
In this work, we adopt the transversely isotropic constitutive model proposed in \cite{guccione1995finite}, known as the \textit{Guccione} relation, with an additional term penalizing large volume variations, so that the passive term of the Piola-Kirchoff stress tensor is given by
\begin{equation}\label{eq:passive_piola}
\mathbf{P}_p(\mathbf{F}) = \frac{\partial}{\partial\mathbf{F}} \left( \mathcal{W}_{Guccione} + \mathcal{W}_{vol} \right)
= \frac{\partial}{\partial\mathbf{F}} \left( \frac{C}{2}(e^Q-1) + \frac{K}{2}(J-1)\text{log}(J) \right),
\end{equation}
with the following form for $Q$ to describe transverse isotropy with respect to the fiber coordinate system (with fibers in the $\mathbf{e}_1$-direction),
\begin{equation*}
Q = b_{f} E_{ff}^2 + b_{s} E_{ss}^2 + b_{n} E_{nn}^2 + b_{fs}(E_{fs}^2 + E_{sf}^2) + b_{fn}(E_{fn}^2 + E_{nf}^2) + b_{sn}(E_{sn}^2 + E_{ns}^2).
\end{equation*}
Here, $E_{ij}$, $i,j\in\{f,s,n\}$ are the components of the Green-Lagrange strain tensor $\mathbf{E}$, the material constant $C>0$ is used for scaling the stresses and the coefficients $b_f$, $b_s$, $b_n$ are related to the material stiffness in the fiber, sheet and cross-fiber directions, respectively. Finally, the bulk modulus $K>0$ is the penalization term controlling the incompressibility of the myocardial tissue. Since active properties are time-dependent and anisotropic \cite{lin1998multiaxial} (with more active stress generated along the local muscle fiber direction), we model the tissue stretch along the reference fiber direction $\mathbf{f}_0\in\mathbb{R}^3$ only and define
\begin{equation}\label{eq:active_piola}
\mathbf{P}_a(\mathbf{F}) = \mathbf{T}_a(\mathbf{X},t;\bm\mu)(\mathbf{Ff}_0\otimes\mathbf{f}_0),
\end{equation}
where $\mathbf{T}_a$ represents the active tension generated at cellular level. A surrogate model for the active tension, introduced to avoid the coupling with the electrophysiology model, is described in Sec.~\ref{sec:Ta}. From now on, since the deformation gradient can be calculated as $\mathbf{F} = \mathbf{I} + \nabla_0\mathbf{u}$, we write $\mathbf{P} = \mathbf{P}(\mathbf{u})$. 

\begin{figure}
	\centering
	\includegraphics[width=0.35\textwidth]{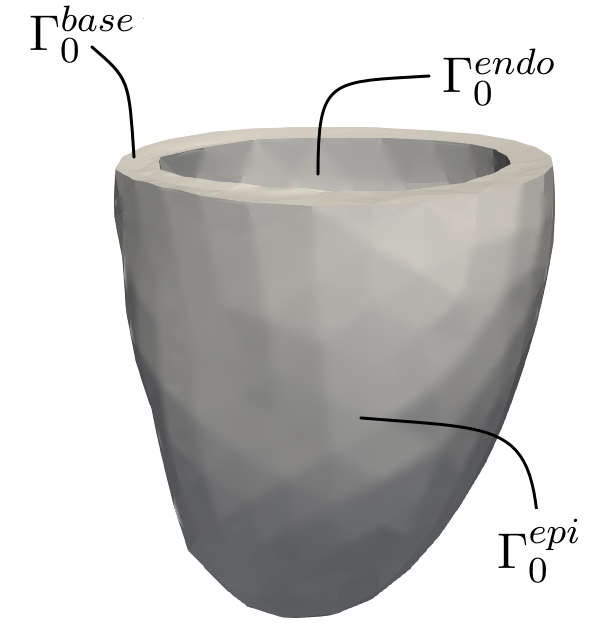}
	\caption{Patient-specific unloaded left ventricle geometry.} 
	\label{fig:LV_BCs}
\end{figure}

The strong formulation of the nonlinear parameterized initial-boundary value problem for cardiac mechanics we consider reads as follows: given $\bm\mu\in\mathcal{P}$, find the displacement field $\mathbf{u}(\bm\mu)\colon\Omega_0\times[0,T)\rightarrow\mathbb{R}^3$ such that
\begin{equation}\label{eq:3D}
\left\{ \begin{array}{lllr}
\rho_0 \ddot{\mathbf{u}}(\bm\mu) - \nabla_0\cdot\mathbf{P}(\mathbf{u}(\bm\mu);\bm\mu) = \mathbf{0} && \text{in } \Omega_0\times(0,T),\\
\mathbf{P}(\mathbf{u}(\bm\mu);\bm\mu){\mathbf{N}} = p_{LV}(t;\bm\mu)\lVert J\mathbf{F}^T(\mathbf{u}(\bm\mu);\bm\mu)\mathbf{N}\rVert\mathbf{v}(t)
&& \text{on } \Gamma_0^{base}\times(0,T), \\
\mathbf{P}(\mathbf{u}(\bm\mu);\bm\mu){\mathbf{N}} + \mathbf{K}^{epi}\mathbf{u}(\bm\mu) + \mathbf{C}^{epi}\dot{\mathbf{u}}(\bm\mu) = 0
&& \text{on } \Gamma_0^{epi}\times(0,T),\\
\mathbf{P}(\mathbf{u}(\bm\mu);\bm\mu){\mathbf{N}} = - p_{LV}(t;\bm\mu)J\mathbf{F}^{-T}(\mathbf{u}(\bm\mu)){\mathbf{N}} && \text{on } \Gamma_0^{endo}\times(0,T),\\
\mathbf{u}(\bm\mu) = \mathbf{u}_0(\bm\mu); ~~\dot{\mathbf{u}}(\bm\mu) = \dot{\mathbf{u}}_0(\bm\mu) && \text{in } \Omega_0\times\{0\},
\end{array} \right.
\end{equation}
where the computational boundary $\partial\Omega_0$ is divided into the inner endocardium $\Gamma_0^{endo}$, the outer epicardium $\Gamma_0^{epi}$ and the ventricular base $\Gamma_0^{base}$, the latter representing the artificial boundary resulting from truncation of the heart below the valves in a short axis plane, see Figure~\ref{fig:LV_BCs}. The boundary conditions on $\Gamma_0^{base}$ are energy-consistent and provide an explicit expression for the stresses at the base \cite{regazzoni2020machine}, being
\begin{equation*}
\mathbf{v}(t) = \frac{ \int_{\Gamma_0^{endo}} J\mathbf{F}^{-T}(\mathbf{u}(\bm\mu)){\mathbf{N}} d\Gamma }{\int_{\Gamma_0^{base}}\lVert J\mathbf{F}^{-T}(\mathbf{u}(\bm\mu)){\mathbf{N}} \rVert d\Gamma}.
\end{equation*} 
The Robin boundary conditions at the epicardium aim at modeling the interaction between the ventricle and the pericardium \cite{gerbi2019monolithic}, that is the fibroelastic sac containing the heart, and are given by $\mathbf{K}^{epi} = K_{\perp}(\mathbf{N}\otimes\mathbf{N}) + K_{\parallel}(\mathbf{I}-\mathbf{N}\otimes\mathbf{N})$ and $\mathbf{C}^{epi} = C_{\perp}(\mathbf{N}\otimes\mathbf{N}) + C_{\parallel}(\mathbf{I}-\mathbf{N}\otimes\mathbf{N})$, where the local values of stiffness $\mathbf{K}^{epi}$ and viscosity $\mathbf{C}^{epi}$ of the epicardial tissue, in the normal ($\perp$) and tangential ($\parallel$) directions, are reported in Table~\ref{tab:reference_values_3D}. Finally, Neumann boundary conditions account for the action of the blood pressure $p_{LV}(t;\bm\mu)$ at the endocardium.

To provide meaningful numerical simulations of the left ventricle activity between two consecutive heartbeats, and then to characterize the complete cardiac cycle from a mechanical point of view, we rely on a suitable, lumped-parameter model for blood circulation, presented in Sec.~\ref{sec:0D}. This strategy allows us to take into account the action of the blood pressure inside the chamber.


\subsection{0D blood external circulation model}\label{sec:0D}

\begin{figure}
	\centering
	\includegraphics[width=0.825\textwidth]{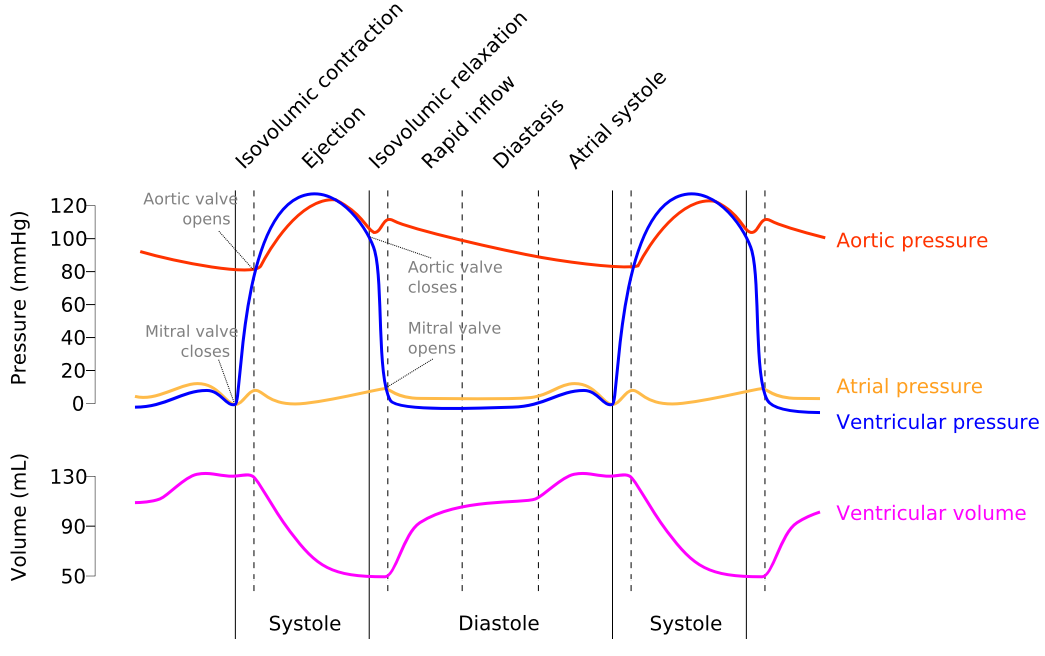}
	\vspace{-0.15cm}
	\caption{Wiggers diagram, adapted from \url{https://commons.wikimedia.org/w/index.php?curid=50317988}, illustrating the events taking place over the cardiac cycle.}
	\label{fig:Wiggers}
\end{figure}

Several hemodynamics models have been proposed in the literature to account for the presence of blood inside the cardiac chamber, see, e.g., \cite{regazzoni2020cardiac, blanco20103d, westerhof2009arterial, nordsletten2011coupling}, just to mention a few examples. Among these, in the context of coupled problems, lumped-parameter fluid models have been extensively considered \cite{hirschvogel2017monolithic, pfaller2020using, mollero2018multifidelity}, since they provide good approximation results at a greatly reduced cost. In this work, we adopt the following 0D model, as done in \cite{gerbi2019monolithic, regazzoni2020machine}. Starting with systole, that is the phase in which the ventricle is full of blood and both the mitral valve and the aortic valve are closed, the four phases of the cardiac cycle can be described as follows (see Wiggers diagram \cite{wiggers1951physiology}, reported in Figure~\ref{fig:Wiggers}):
\begin{enumerate}
	\item\label{it:iso_contraction} \textit{isovolumetric contraction:} the endocardial pressure rapidly grows from the end-diastolic pressure $p_{ED}$ to the value measured in the aorta, in such a way that the volume remains unchanged;
	\item\label{it:ejection} \textit{ejection:} as soon as the aortic valve opens, the ejection phase starts and the evolution of the pressure $p_{LV}(t;\bm\mu)$ is governed by a two-element windkessel model \cite{westerhof2009arterial}, with capacitance $C_p$ and resistance $R_p$:
	\begin{equation}\label{eq:windkessel}
	\left\{
	\begin{array}{ll}
	C_p~\dot{p}_{LV}(t;\bm\mu) = -\dfrac{p_{LV}(t;\bm\mu)}{R_p} - \dot{V}_{LV}(\mathbf{u}(\bm\mu),t;\bm\mu), & t\in(T_{AVO},T_{AVC}], \\
	p(T_{AVO};\bm\mu) = p_{AVO}.
	\end{array}
	\right.
	\end{equation}
	Here, $T_{AVO}$ and $T_{AVC}$ are the aortic valve opening and closing times, respectively, and $p_{AVO}$ is the pressure measured in the aorta at the beginning of the ejection phase. This phase is characterized by a decrement of the volume due to the contraction of the ventricle;
	\item\label{it:iso_relaxation} \textit{isovolumetric relaxation:} when the aortic valve closes, the ventricle relaxes and the pressure drops. As both the ventricular valves are closed, no change of volume is experienced;
	\item\label{it:filling} \textit{filling:} finally, as the pressure inside the ventricle falls below that in the atrium, the mitral valve opens and the ventricle begins to fill again, so that the pressure linearly increases to the end-diastolic pressure $p_{ED}$, concluding the cardiac cycle.	
\end{enumerate}

\subsection{The full order model}

The mechanics and the blood circulation core models described so far, which mutually exchange pressure and volume, must be suitable coupled to provide physically meaningful simulations of the cardiac cycle, leading to a 3D-0D coupled structure-windkessel problem. In this section, we outline the corresponding full order model obtained by relying on the FEM in space and on implicit time schemes, which avoid restrictions on the time step due to the highly nonlinear terms of the strain energy density function.

Let $\mathcal{T}_h$ be an hexahedral mesh on the reference domain $\Omega_0\subset\mathbb{R}^3$ and $\mathbb{Q}_r(\tau)$ be the set of polynomials of degree smaller than or equal to $r\geq1$ over a mesh element $\tau\in\mathcal{T}_h$. Given the finite-dimensional space of real-valued functions
\begin{equation*}
\mathcal{X}_h^r = \left\{v\in C^0(\bar{\Omega}_0)~\colon v\lvert_{\tau}\in\mathbb{Q}_r(\tau)~~\forall\tau\in\mathcal{T}_h\right\},
\end{equation*}
we define the finite element (FE) space of degree $r\geq1$ as \vspace{-0.1cm}
\begin{equation*}
V_h = \left\{\bm\eta\in [H^1(\Omega_0)]^3~|~\bm\eta=\mathbf{0}~\text{on}~\Gamma_0^D\right\} \cap [\mathcal{X}_h^r]^3, \vspace{-0.1cm}
\end{equation*}
whose dimension $N_h=3\dim(\mathcal{X}_h^r)$ corresponds to the total number of structural dofs. Furthermore, we introduce a uniform partition $\{t^0,\dots,t^{N_t}\}$ of the time interval $(0,T)$, corresponding to the duration of a single heartbeat, with time step $\Delta t$. The vector of nodal displacements of the fully-discretized problem and the left ventricular pressure at time $t^n$, for $n=1,\dots,N_t$, are denoted as $\mathbf{u}_h^n\in\mathbb{R}^{N_h}$ and $p_{LV}^{n}(\bm\mu)$, respectively. Time derivatives computed at time $t^n$ are approximated as 
\begin{align*}
\partial_t\mathbf{u}(t^n) \approx \frac{\mathbf{u}_h^{n}-\mathbf{u}_h^{n-1}}{\Delta t}, &&
\partial_t^2\mathbf{u}(t^n) \approx \frac{\mathbf{u}_h^{n}-2\mathbf{u}_h^{n-1}+\mathbf{u}_h^{n-2}}{\Delta t^2}.
\end{align*}
To ease the notation of what follows, we define \vspace{-0.1cm}
\begin{equation*}
(\star_h^{n,(k)}) := (\mathbf{u}_h^{n,(k)}(\bm\mu),p_{LV}^{n,(k)}(\bm\mu),t^{n};\bm\mu), \vspace{-0.1cm}
\end{equation*}
where the superscript $(k)$ denotes quantities computed at the $k$-th iteration of Newton method used for the solution of the algebraic nonlinear system arising at each time step. Note that $p_{LV}^{n,(k)}(\bm\mu):=p_{LV}^n(\bm\mu)$, $\forall k\geq0$, during the non-isochoric phases.

During ventricular ejection (phase \ref{it:ejection}) and filling (phase \ref{it:filling}), the structural and the circulation problems are segregated, meaning that the two models are solved one after the other. In particular, in the ejection phase, the current pressure $p_{LV}^n(\bm\mu)$ is updated by solving the two-element windkessel model (\ref{eq:windkessel}) before addressing the mechanics problem. For simplicity, we assume 
\begin{equation*}
\dot{V}_{LV}(\mathbf{u}_h^n(\bm\mu),t^n;\bm\mu)\approx\frac{V_{LV}(\mathbf{u}_h^{n-1}(\bm\mu),t^{n-1};\bm\mu)-V_{LV}(\mathbf{u}_h^{n-2}(\bm\mu),t^{n-2};\bm\mu)}{\Delta t}, \vspace{-0.1cm}
\end{equation*} 
where the ventricular volume at time $t^j$, for $j\in\{n-1,n-2\}$, is computed as  \vspace{-0.1cm}
\begin{equation}\label{eq:volumeVentr}
V_{LV}(\mathbf{u}_h^j(\bm\mu),t^j;\bm\mu) = \frac{1}{3}\int_{\Gamma_0^{endo}} J(\mathbf{X}+\mathbf{u}_h^j(\bm\mu)-\mathbf{b}_h^j(\bm\mu))\cdot\mathbf{F}^{-T}\mathbf{N}d\Gamma_0,
\end{equation}
being $\mathbf{b}_h^j(\bm\mu) = \frac{1}{\lvert\Gamma_0^{base}\rvert}\int_{\Gamma_0^{base}} (\mathbf{X}+\mathbf{u}_h^j(\bm\mu))d\Gamma_0$; for further details on the derivation of  formula \eqref{eq:volumeVentr} we refer to \cite{regazzoni2020machine}. The corresponding problem at time $t^n$ for the unknown $\mathbf{u}_h^{n}(\bm\mu)$, for $n=1,\dots,N_t$, is given by the nonlinear system
\begin{equation*}\label{eq:ejection_filling_problem}
\mathbf{R}(\mathbf{u}_h^{n}(\bm\mu),p_{LV}^n(\bm\mu),t^{n};\bm\mu) = \mathbf{0}\quad \text{in } \mathbb{R}^{N_h}
\end{equation*}
and is solved by means of the Newton method, leading to a sequence of linear systems of the form
\begin{equation*}
\partial_\mathbf{u}\mathbf{R}(\star_h^{n,(k)})\delta\mathbf{u}_h^{(k)}(\bm\mu) = - \mathbf{R}(\star_h^{n,(k)}),\quad k\geq0,
\end{equation*}
where $\partial_\mathbf{u}\mathbf{R}$ is the directional derivative of the structural residual. At each iteration $k$, the current solution is thus updated as $\mathbf{u}_h^{n,(k+1)}(\bm\mu) = \mathbf{u}_h^{n,(k)}(\bm\mu) + \delta\mathbf{u}_h^{(k)}(\bm\mu)$.

On the other hand, during isovolumetric contraction (phase \ref{it:iso_contraction}) and isovolumetric relaxation (phase \ref{it:iso_relaxation}), the elastodynamics problem is solved together with the volume constraint $V_{LV}^{n} = V_{LV}^{n-1}$. This results in a nonlinear saddle-point system for the unknowns $\mathbf{u}_h^n(\bm\mu)$ and $p_{LV}^n(\bm\mu)$ of the form \vspace{-0.1cm}
\begin{equation*}
\left\{
\begin{array}{l}
\mathbf{R}(\mathbf{u}_h^n(\bm\mu),p_{LV}^n(\bm\mu),t^n;\bm\mu) = \mathbf{0},\\
V_{LV}(\mathbf{u}_h^n(\bm\mu),t^n;\bm\mu) = V_{LV}(\mathbf{u}_h^{n-1}(\bm\mu),t^{n-1};\bm\mu),
\end{array}
\right.
\end{equation*}
that can be solved by means of the Schur complement reduction \cite{benzi2005numerical}. By applying Newton method, we end up with the following linear system \vspace{-0.1cm}
\begin{equation*}
\left( \begin{array}{l}
~\partial_\mathbf{u}\mathbf{R}(\star_h^{n,(k)}) \quad \partial_p\mathbf{R}(\star_h^{n,(k)}) \\ 
\partial_\mathbf{u}\text{R}^{vol}(\star_h^{n,(k)}) \qquad 0 \end{array} \right)
\left( \begin{array}{c}
\delta\mathbf{u}_h^{(k)}(\bm\mu)\\
\delta p_{LV}^{(k)}(\bm\mu)
\end{array} \right) = -
\left( \begin{array}{c}
\mathbf{R}(\star_h^{n,(k)}) \\
\text{R}^{vol}(\star_h^{n,(k)})
\end{array} \right) \vspace{-0.1cm}
\end{equation*}
at each iteration $k\geq0$, where $\text{R}^{vol}\in\mathbb{R}$ is the residual related to the volume constraint.

To summarize, the discrete nonlinear parameterized FOM for the coupled problem can be written as: given $\bm\mu\in\mathcal{P}$, for $n=1,\dots,N_t$, find $\mathbf{u}_h^n(\bm\mu)\in\mathbb{R}^{N_h}$ and $p_{LV}^{n}(\bm\mu)>0$ such that
\begin{equation}\label{eq:FOM_3D-0D}
\left[
\begin{array}{c}
\mathbf{R}(\mathbf{u}_h^n(\bm\mu),p_{LV}^n(\bm\mu),t^n;\bm\mu)\\
\text{R}^{vol}(\mathbf{u}_h^{n}(\bm\mu),p_{LV}^{n}(\bm\mu),t^{n};\bm\mu)
\end{array}
\right] = \mathbf{0},
\end{equation}
where $\text{R}^{vol}$ is discarded during phases \ref{it:ejection} and \ref{it:filling} of the cardiac cycle.
 
The numerical solution of problem (\ref{eq:FOM_3D-0D}) entails huge computational costs as soon as $N_h$ (depending on the computational mesh and on the discretization scheme) becomes too large. This is extremely challenging, if not prohibitive, when the repeated solution to the forward problem is required, such as in the context of UQ, parameter estimation or model calibration.


\section{Deep-HyROMnet for the cardiac mechanics problem} \label{sec:LV-Deep-HyROMnet}

With the aim of reducing the computational burden associated with the FOM, we address the solution to the 3D-0D coupled problem described in Sec.~\ref{sec:3D-0D} by means of the \textit{deep hyper-reduced order model network} (Deep-HyROMnet) technique \cite{cicci2021DeepHyROMnet} for the efficient solution to time-dependent, nonlinear parameterized PDEs; further details on this reduction strategy are provided in Appendix~\ref{appendix:DeepHyROMnet}. We point out that, since blood circulation is modeled through a lumped-parameter model and the volume constraint implies only a few additional dofs to the mechanics problem, the 0D circulation model does not need to be reduced, similarly to the approach adopted in \cite{pfaller2020using}.
 
Based on the Galerkin-RB method \cite{quarteroni2016reduced, hesthaven2016certified}, we aim at approximating the elements of the high-fidelity discrete solution manifold
\begin{equation*}
\mathcal{M}_h = \{\mathbf{u}_h^n(\bm\mu)\in\mathbb{R}^{N_h}, n=1,\dots,N_t ~|~ \bm\mu\in\mathcal{P} \}
\end{equation*}
by means of a linear combination of (possibly few) global, problem-dependent, basis functions. For $n=1,\dots,N_t$, the reduced displacement $\mathbf{u}_N^n(\bm\mu)\in\mathbb{R}^{N}$ such that
\begin{equation*}
\mathbf{Vu}_N^n(\bm\mu) \approx \mathbf{u}_h^n(\bm\mu),
\end{equation*}
where $\mathbf{V}\in\mathbb{R}^{N_h\times N}$ ($N\ll N_h$) is the matrix collecting column-wise the nodal values of the RB functions, is found by solving a low-dimensional nonlinear problem obtained by requiring the fulfillment of a suitable orthogonality criterion. In this work, the reduced basis $\mathbf{V}$ is built by performing POD (see Appendix~\ref{appendix:POD}) on the snapshots matrix of mechanical displacements, i.e.
\begin{equation*}
\mathbf{S}_u = \left[\mathbf{u}_h^1(\bm\mu_1)~|~\dots~|~\mathbf{u}_h^{N_t}(\bm\mu_1)~|~\dots~|~\mathbf{u}_h^1(\bm\mu_{n_s})~|~\dots~|~\mathbf{u}_h^{N_t}(\bm\mu_{n_s})\right],
\end{equation*}
for randomly sampled parameter values $\bm\mu_1,\dots,\bm\mu_{n_s}$.  Performing a Galerkin projection of the residual (\ref{eq:FOM_3D-0D})$_1$ of the full-order structural model  onto the reduced subspace spanned by the columns of $\mathbf{V}$, we obtain the low-dimensional problem 
\begin{equation}\label{eq:ROM_3D-0D}
\left[
\begin{array}{c}
\mathbf{V}^T\mathbf{R}(\mathbf{Vu}_N^n(\bm\mu),p_{LV}^n(\bm\mu),t^n;\bm\mu)\\
\text{R}^{vol}(\mathbf{Vu}_N^{n}(\bm\mu),p_{LV}^n(\bm\mu),t^n;\bm\mu)
\end{array}
\right] = \mathbf{0}.
\end{equation}
As done before, to ease the notation, we define
\begin{equation*}
(\star_N^{n,(k)}) := (\mathbf{Vu}_N^{n,(k)}(\bm\mu),p_{LV}^{n,(k)}(\bm\mu),t^{n};\bm\mu).
\end{equation*}
The corresponding reduced Newton system at time $t^n$, for $n=1,\dots,N_t$, thus reads:
\begin{itemize}	
	\item for the ejection and filling phases: given an initial guess $\mathbf{u}_N^{n,(0)}(\bm\mu)$, find $\mathbf{u}_N^{n,(k)}(\bm\mu)$ such that, for $k\geq0$, 
	\begin{equation*}
	\left\{
	\begin{array}{l}
	\mathbf{V}^T \partial_\mathbf{u}\mathbf{R}(\star_N^{n,(k)})\mathbf{V}\delta \mathbf{u}_N^{(k)} = - \mathbf{V}^T\mathbf{R}(\star_N^{n,(k)}),\\
	\mathbf{u}_N^{n,(k+1)}(\bm\mu) = \mathbf{u}_N^{n,(k)}(\bm\mu) + \delta\mathbf{u}_N^{(k)}(\bm\mu),
	\end{array} 
	\right.
	\end{equation*}
	until $\lVert \mathbf{V}^T\mathbf{R}(\star_N^{n,(k+1)})\rVert_2 / \lVert\mathbf{V}^T\mathbf{R}(\star_N^{n,(0)})\rVert_2<\varepsilon_{Nwt}$, where $\varepsilon_{Nwt}>0$ is a prescribed tolerance;
	\item for the isovolumetric phases: given initial guesses $\mathbf{u}_N^{n,(0)}(\bm\mu)$ and $p_{LV}^{n,(0)}(\bm\mu)$, find $\mathbf{u}_N^{n,(k)}(\bm\mu)$ and $p_{LV}^{n,(k)}(\bm\mu)$ such that, for $k\geq0$, 
	\begin{equation*}
	\left( \begin{array}{l}
	\mathbf{V}^T\partial_\mathbf{u}\mathbf{R}(\star_N^{n,(k)})\mathbf{V} \quad \mathbf{V}^T\partial_p\mathbf{R}(\star_N^{n,(k)})\\ 
	~~\partial_\mathbf{u}\text{R}^{vol}(\star_N^{n,(k)})\mathbf{V} \quad \hspace{8mm} 0 
	\end{array} \right)
	\left( \begin{array}{c}
	\delta\mathbf{u}_N^{(k)}(\bm\mu)\\
	\delta p_{LV}^{(k)}(\bm\mu)
	\end{array} \right) = -
	\left( \begin{array}{c}
	\mathbf{V}^T\mathbf{R}(\star_N^{n,(k)}) \\
	\text{R}^{vol}(\star_N^{n,(k)})
	\end{array} \right),
	\end{equation*}
	then update 
	\[
	\mathbf{u}_N^{n,(k+1)}(\bm\mu) = \mathbf{u}_N^{n,(k)}(\bm\mu) + \delta\mathbf{u}_N^{(k)}(\bm\mu)
	\]
	 and 
	 \[
	 p_{LV}^{n,(k+1)}(\bm\mu) = p_{LV}^{n,(k)}(\bm\mu) + \delta p_{LV}^{(k)}(\bm\mu),
	 \] until $\lVert \mathbf{V}^T\mathbf{R}(\star_N^{n,(k+1)})\rVert_2 / \lVert\mathbf{V}^T\mathbf{R}(\star_N^{n,(0)})\rVert_2<\varepsilon_{Nwt}$, where $\varepsilon_{Nwt}>0$ is a prescribed tolerance.
\end{itemize}
As initial guess we choose $\mathbf{u}_N^{0,(0)}(\bm\mu)=\mathbf{u}_0(\bm\mu)$, given by the initial condition (\ref{eq:3D})$_5$, and $\mathbf{u}_N^{n,(0)}(\bm\mu) = \mathbf{u}_N^{n-1}(\bm\mu)$, for $n=1,\dots,N_t$.

Since the reduced arrays $\mathbf{V}^T\mathbf{R}\in\mathbb{R}^{N\times1}$, $\mathbf{V}^T\partial_\mathbf{u}\mathbf{R}\mathbf{V}\in\mathbb{R}^{N\times N}$, $\mathbf{V}^T\partial_p\mathbf{R}\in\mathbb{R}^{N\times1}$, $\text{R}^{vol}\in\mathbb{R}^{1\times1}$ and $\partial_\mathbf{u}\text{R}^{vol}\mathbf{V}\in\mathbb{R}^{1\times N}$ are evaluated on the current solutions $\mathbf{Vu}_N^{n,(k)}(\bm\mu)$ and $p_{LV}^{n,(k)}(\bm\mu)$, they have to be com\-put\-ed for every new $k\geq0$ and $n=1,\dots,N_t$. However, due to nonlinearity, the corresponding high-fidelity arrays must be assembled at each Newton iteration before projecting them onto the reduced subspace, thus entailing a computational cost that still depends on $N_h$. To overcome this limitation, suitable hyper-reduction techniques should be taken into account, in order to provide approximations of the nonlinear terms that are independent of the FOM dimension.

The discrete empirical interpolation method (DEIM) \cite{chaturantabut2010nonlinear}, represents a standard hyper-reduction technique very often used in a POD-Galerkin setting. It depends on the assembling of the nonlinear quantities onto a reduced mesh obtained as a subset of the original one. Nonetheless, when applied in the context of cardiac mechanics, this strategy still suffers from severe computational burdens, as a large reduced mesh is required to correctly capture the great variability of the (nonlinear) residual vectors \cite{cicci2021cardiacDEIM}.

With the aim of avoiding the assembling stage and thus overcome the computational bottleneck associated with DEIM, we perform a deep learning-based approximation of the reduced nonlinear terms. Given the triplets 
\begin{equation*}
\bm\vartheta = (\bm\mu,t^n,k)\in\mathcal{P}\times\{t^1,\dots,t^{N_t}\}\times\mathbb{N}^+
\end{equation*}
made of the input parameters $\bm\mu\in\mathcal{P}$, the current time step $t^n$ and the Newton iteration $k\geq0$, we efficiently compute the $N$-dimensional ROM operators evaluated on $(\star_N^{n,(k)})$ by exploiting the DNN architecture described in the Appendix~\ref{appendix:DeepHyROMnet} to learn the following nonlinear maps:
\begin{align*}
\bm\rho_N\colon(\bm\mu,t^n,k) & \longmapsto 
\bm\rho_N(\bm\mu,t^n,k) \approx\mathbf{V}^T\mathbf{R}(\star_N^{n,(k)}),\\
\bm\iota_N\colon(\bm\mu,t^n,k) & \longmapsto \bm\iota_N(\bm\mu,t^n,k) \approx\mathbf{V}^T\partial_\mathbf{u}\mathbf{R}(\star_N^{n,(k)})\mathbf{V},\\
\bm\pi_N\colon(\bm\mu,t^n,k) & \longmapsto \bm\pi_N(\bm\mu,t^n,k) \approx\mathbf{V}^T\partial_p\mathbf{R}(\star_N^{n,(k)}),\\
\bm\upsilon_N\colon(\bm\mu,t^n,k) & \longmapsto [\bm\upsilon_N^{\partial R}(\bm\mu,t^n,k), \bm\upsilon_N^R(\bm\mu,t^n,k)] \approx [\partial_\mathbf{u}\text{R}^{vol}(\star_N^{n,(k)})\mathbf{V},\text{R}^{vol}(\star_N^{n,(k)})].
\end{align*}
This procedure guarantees an efficient decomposition into a costly (offline) training phase, which is performed once and for all, and an inexpensive (online) testing phase, during which the problem solution is computed for a specific input vector $\bm\mu\in\mathcal{P}$. During the offline phase, we need to collect FOM snapshots for the construction of the reduced basis $\mathbf{V}$. Then, reduced nonlinear data are collected by performing ROM simulations for a new set of input parameter values, i.e. different from the ones used for basis construction, and the neural networks (NNs) are trained. Online, for each new instance of the input parameter, the output of the NNs is evaluated in order to assemble the reduced Newton system, thus recovering the efficiency of the reduced model.

More precisely, the online stage reads as follows: given $\bm\mu\in\mathcal{P}$, for $n=1,\dots,N_t$, given $\mathbf{u}_N^{n,(0)}(\bm\mu)$ and $p_{LV}^{n,(0)}(\bm\mu)$, for $k\geq0$, find $\mathbf{\delta u}_N^{(k)}(\bm\mu)\in\mathbb{R}^{N}$ and $\delta p_{LV}^{(k)}>0$ such that
\begin{equation}\label{eq:DeepLin1}
\bm\iota_N(\bm\mu,t^n,k)\delta \mathbf{u}_N^{(k)} = - \bm\rho_N(\bm\mu,t^n,k)
\end{equation}
for the ejection and filling phase, and by 
\begin{equation}\label{eq:DeepLin2}
\left( \begin{array}{l}
~\bm\iota_N(\bm\mu,t^n,k) \quad \bm\pi_N(\bm\mu,t^n,k)\\ 
\bm\upsilon_N^{\partial R}(\bm\mu,t^n,k) \hspace{13mm} 0
\end{array} \right)
\left( \begin{array}{c}
\delta\mathbf{u}_N^{(k)}\\
\delta p_{LV}^{(k)}
\end{array} \right) = -
\left( \begin{array}{c}
\bm\rho_N(\bm\mu,t^n,k) \\
\bm\upsilon_N^R(\bm\mu,t^n,k)
\end{array} \right)
\end{equation}
for the isovolumetric phases, until $\lVert \bm\rho_N(\bm\mu,t^n,k+1)\rVert_2/\lVert \bm\rho_N(\bm\mu,t^n,0)\rVert_2<\varepsilon_{Nwt}$.

Thanks to the employed DNN architectures, the linear systems (\ref{eq:DeepLin1}) and (\ref{eq:DeepLin2}) are assembled in an extremely efficient way, i.e. requiring $\mathcal{O}(10^{-3})$~s, while both FOM and POD-Galerkin ROMs require $\mathcal{O}(10^{-1})$~s or even $\mathcal{O}(1)$~s for each Newton iteration. Since this operation is performed about $N_tN_{nwt}$ times during each cardiac cycle (where $N_{nwt}$ is the average number of Newton iterations per time step), relying on the Deep-HyROMnet strategy allows us to achieve speed-ups with respect to the FOM of more than two orders of magnitude regarding CPU time, as shown in Sec.~\ref{sec:numerical_results}. Traditional hyper-reduction techniques as DEIM would not allow such a computational gain.


\section{Numerical results}\label{sec:numerical_results}

In this section we present the numerical results obtained using our Deep-HyROMnet strategy for the solution to the 3D-0D structure-windkessel model in both physiological and pathological scenarios. Regarding the high-fidelity model, we point out that quadratic ($\mathbb{Q}_2$) FE are commonly used \cite{land2015verification} when dealing with cardiac mechanics, especially in a nearly-incompressible regime, due to possible instabilities. However, due to the huge computational costs entailed, we rely on trilinear ($\mathbb{Q}_1$) FE, which proved to be sufficiently accurate for the purposes at hand and less expensive, despite considering suitable refined meshes. We point out that no instabilities have been observed. Moreover, we recall that the reduction strategy, acting at the algebraic level, works irrespectively of the chosen FE degree. 

\begin{figure}[b]
	\centering
	\includegraphics[width=\textwidth]{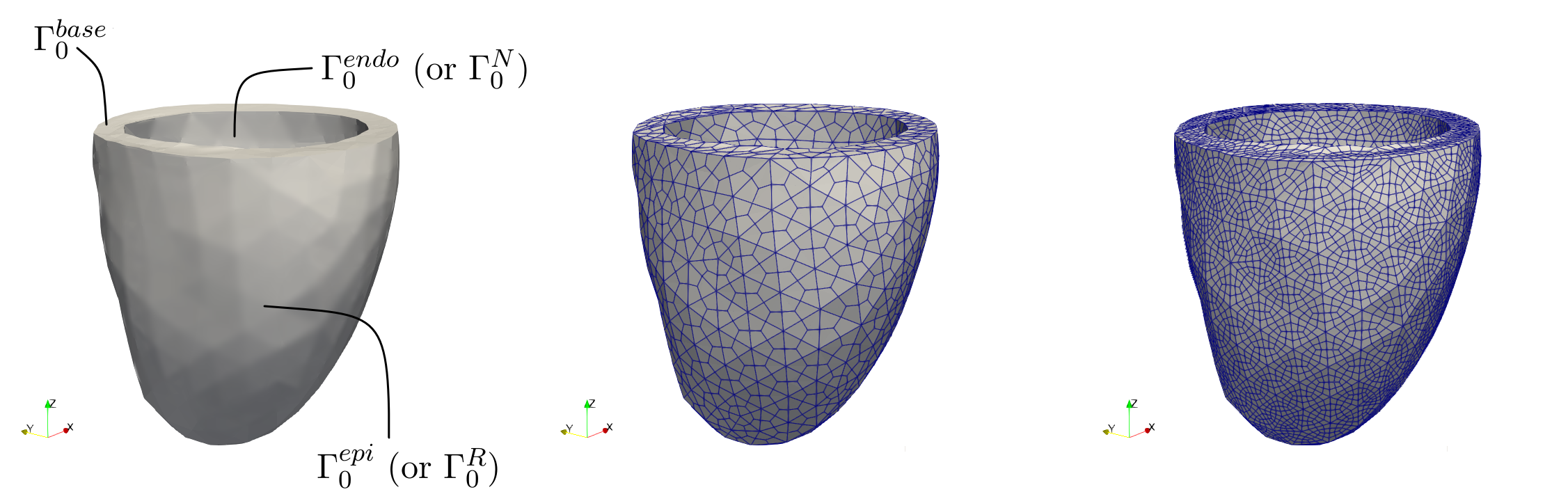}
	\caption{Patient-specific loaded left ventricle geometry (left) and computational grids (center and right).} 
	\label{fig:LV_loaded_BCs}
\end{figure}

In Figure~\ref{fig:LV_loaded_BCs} we report the computational geometry obtained when the ventricle is loaded by a value of pressure corresponding to the end diastolic pressure, in our case $p_{ED}=15$~mmHg, and the hexahedral meshes used.

\begin{remark}
	To correctly start the numerical simulation, we need to find the end-diastolic con\-fig\-uration of the left ventricle and to use the corresponding displacement as initial condition for our problem. This is done by solving the quasi-static problem (obtained from (\ref{eq:3D}) by setting to zero the time dependent terms, see, e.g., \cite{regazzoni2020cardiac}) on the reference configuration, so that the resulting initial displacement depends on the input parameters. For practical reasons, we solve the initial displacement problem once and for all given the reference values of the input parameters, reported in Tables~\ref{tab:reference_values_3D} and \ref{tab:reference_values_0D}, that we collectively denote as $\tilde{\bm\mu}$, so that the initial conditions are  $\mathbf{u}_{h,0} = \mathbf{u}_{h,0}(\tilde{\bm\mu})$ and $\dot{\mathbf{u}}_{h,0} = \dot{\mathbf{u}}_{h,0}(\tilde{\bm\mu})$ for every instance of the parameter vector, both during training (offline stage) and testing (online stage). Nonetheless, a reduced model for the quasi-static problem can be developed in order to take into account different initial conditions as well.
\end{remark}

\begin{figure}
	\centering
	\includegraphics[width=0.65\textwidth]{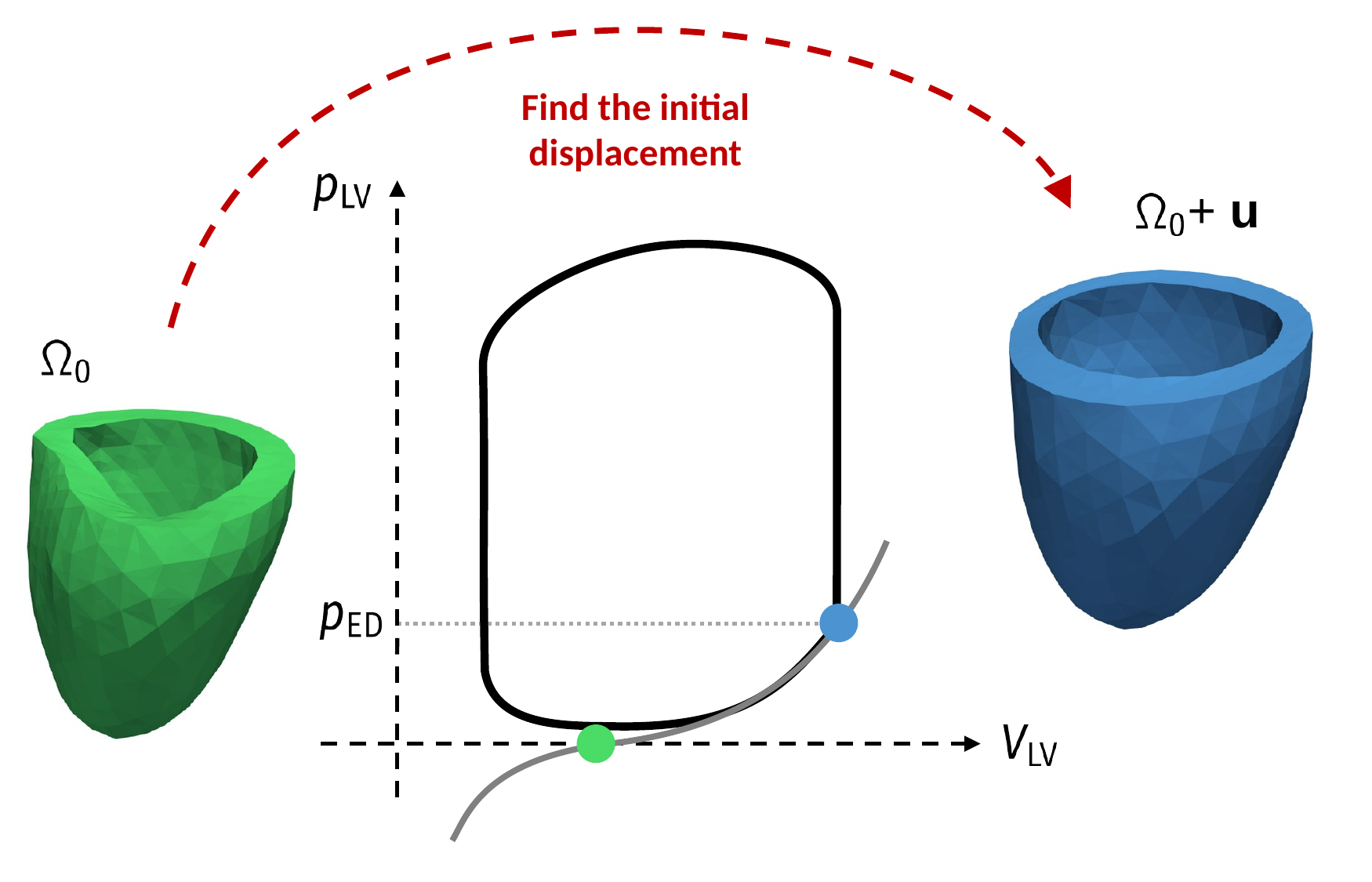}
	\caption{Sketch of the strategy to initialize the simulation, adapted from \cite{regazzoni2020cardiac}. The black line represents the pressure-volume loop, while the gray line is the Klotz curve \cite{klotz2006single}.}
	\label{fig:refConf_initDispl}
\end{figure}

In order to evaluate the accuracy of the ROM with respect to the FOM, the following time-averaged $L^2$-errors of the displacement vector are used
\begin{align}
\epsilon_{abs}(\bm\mu) &= \frac{1}{N_t}\sum_{n=1}^{N_t} \lVert \mathbf{u}_h(\cdot,p^n,t^n;\bm\mu) - \mathbf{Vu}_N(\cdot,p^n,t^n;\bm\mu)\rVert_{2} \label{err_abs}\\
\epsilon_{rel}(\bm\mu) &= \frac{1}{N_t}\sum_{n=1}^{N_t} \frac{\lVert \mathbf{u}_h(\cdot,p^n,t^n;\bm\mu) - \mathbf{Vu}_N(\cdot,p^n,t^n;\bm\mu)\rVert_{2}}{\lVert \mathbf{u}_h(\cdot,p^n,t^n;\bm\mu)\rVert_{2}}, \label{err_rel}
\end{align}
whilst model efficiency is assessed through the CPU time ratio, corresponding to the speed-up achieved by the ROM with respect to the FOM. All the computations have been performed on a PC desktop computer with 3.70GHz Intel Core i5-9600K CPU and 16GB RAM using the code implemented in Python in our software package \texttt{pyfe}$^\text{x}$, which contains a Python binding with the in-house Finite Element library \texttt{life$^\texttt{x}$} (\url{https://lifex.gitlab.io/lifex}), a high-performance C++ library developed within the iHEART project\footnote{iHEART - An Integrated Heart Model for the simulation of the cardiac function, European Research Council (ERC) grant agreement No 740132, P.I. Prof. A. Quarteroni} and based on the \texttt{deal.II} (\url{https://www.dealii.org}) Finite Element core \cite{dealII92}.


\subsection{Parametric setting}\label{sec:Ta}

For all the numerical examples reported in the following sections, we consider a uniform time step $\Delta t = 2.5\cdot10^{-3}$~s for time discretization and set the final time equal to $T=0.8$~s to model a single heartbeat. For the construction of the fiber distribution, we employ the so-called Bayer-Blake-Plank-Trayanova algorithm proposed in \cite{bayer2012novel}, depending on angles coefficients $\bm\alpha^{epi}$, $\bm\alpha^{endo}$, $\bm\beta^{epi}$ and $\bm\beta^{endo}$.

In order to surrogate the input provided to tissue mechanics by the active force generation model, we consider a uniform activation of the cardiac myocytes in the healthy tissue. To be more specific, let
\begin{equation*}
B_r=B(\mathbf{X}_c,r) = \left\{\mathbf{X}\in\Omega_0 ~|~ \lVert \mathbf{X}_c - \mathbf{X} \rVert_2 < r \right\}\subset\mathbb{R}^3
\end{equation*}
be an idealized ischemic region with (fixed) center $\mathbf{X}_c\in\Omega_0$ and radius $r\geq0$, where $r=0$ corresponds to a physiological scenario. The active tension in (\ref{eq:active_piola}) is defined as
\begin{equation*}
\mathbf{T}_a(\mathbf{X},t;\bm\mu) = T_a(t;\bm\mu)\rchi_{B^c_r}(\mathbf{X}),
\end{equation*}
where we assume zero activation in the dofs belonging to the affected region to model the fact that cardiomyocytes inside the necrosis behave as passive conductors. Here $T_a(t;\bm\mu)$ is a prescribed time-dependent function, that can be computed as follows:
\begin{enumerate}
	\item for a fixed set of physiological parameters (see Tables~\ref{tab:reference_values_3D} and \ref{tab:reference_values_0D}), solve the 3D electromechanics (EM) problem coupled with circulation model and an active force generation model for a single heartbeat in the time interval $(0,0.8)$~s. Here, we rely on the model implemented in \cite{regazzoni2020cardiac2}, where the different core models are discretized in space and time by means of the FEM and suitable explicit-implicit differentiation schemes, respectively, and thus solved sequentially (being the time step used for the solution to electrophysiology smaller than that used for the mechanics);		\vspace{-0.1cm}
	\item compute the space-average of the active tension coming from the EM simulation,
	\begin{equation}\label{eq:tau_a^M}
	\tau_a^{EM}(t) = \underset{\mathbf{X}\in\Omega_0}{\text{avg}} T_a^{EM}(\mathbf{X},t),
	\end{equation}	
	and perform a cubic spline interpolation of $\tau_a^{EM}(t)$ to obtain the corresponding time-dependent function $\tau_a^M(t)$, reported in Figure~\ref{fig:LV_Ta}; \vspace{-0.1cm}
	\begin{figure}[b!]
		\centering
		\vspace{-0.25cm}
		\includegraphics[width=0.5\textwidth]{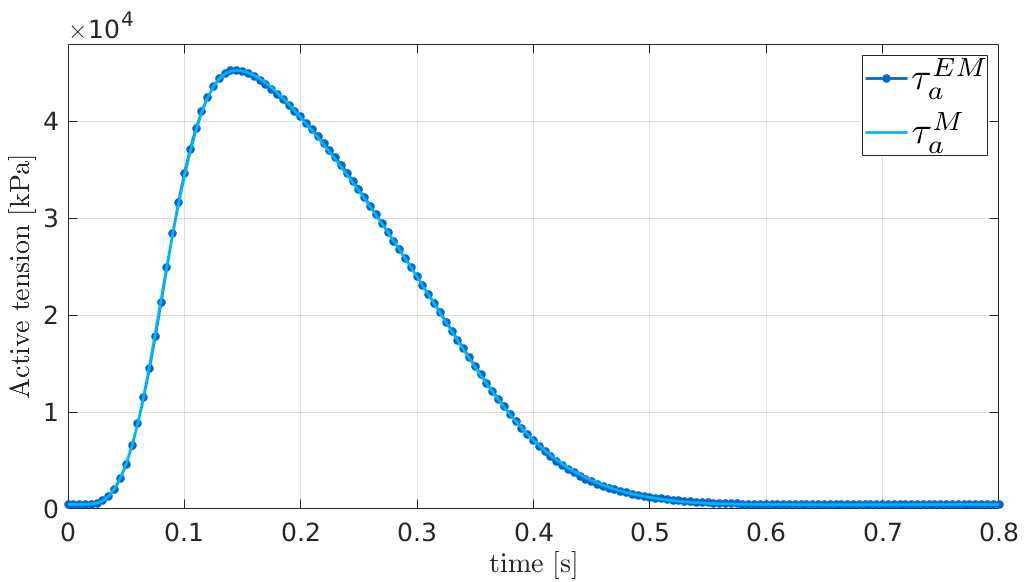}
		\vspace{-0.25cm}
		\caption{Space-averaged active tension computed during a EM simulation and the corresponding fitted curve $\tau_a^{M}(t)$.}
		\label{fig:LV_Ta}
	\end{figure}
	\item finally, for an input parameter $\widetilde{T}_a>0$, define the active tension as
	\begin{equation*}
	T_a(t;\bm\mu) = \frac{\widetilde{T}_a}{\underset{t\in(0,0.8)s}{\max}~\tau_a^{M}(t)}~\tau_a^{M}\hspace{-1mm}\left(0.8~\left(\frac{t}{T} - \floor*{\frac{t}{T}}\right)\right)
	\end{equation*}
	to  account for  parameter-dependence and to model different periodic functions.
\end{enumerate}
To summarize, we define the active tension $T_a(t;\bm\mu)$ used to model the contraction of cardiomyofiber bundles in reference fiber direction by performing a cubic spline interpolation of the average tension computed from the solution to a 3D-0D EM problem and introducing a scaling factor depending on the parameter $\widetilde{T}_a$.


\subsection{Physiological scenario} \label{sec:physio}

First of all, we present the results regarding the Deep-HyROMnet approximation of the FOM solution on physiological, yet challenging, scenarios in cardiac mechanics. In this case, we choose as unknown parameters
\begin{itemize}
	\item the bulk modulus in the passive material law $K\in[2.5\cdot10^{4},7.5\cdot10^{4}]$~Pa,
	\item the resistance of the windkessel model $R_p\in[2.5\cdot10^{7},4.5\cdot10^{7}]~\text{Pa}\cdot\text{s}\cdot\text{m}^{-3}$,
	\item the active tension parameter $\widetilde{T}_a\in[4.5\cdot10^{4},6\cdot10^{4}]$~Pa,
\end{itemize}
that is we set 
$\bm\mu = [K,R_p,\widetilde{T}_a]\in\mathcal{P}\subset\mathbb{R}^3$. 
The parameter space $\mathcal{P}$ is chosen in order to test the accuracy of Deep-HyROMnet in a wide range of scenarios. All other parameters are fixed to their reference values reported in Tables~\ref{tab:reference_values_3D} and \ref{tab:reference_values_0D} in the Appendix.

For the sake of testing, a FOM is built on an hexahedral mesh with $4588$ elements and $6167$ vertices, featuring a high-fidelity dimension equal to $N_h=18501$. During the offline stage, we collect the solution snapshots for 20 parameter samples, each requiring almost 30 minutes using the FOM, and apply POD for the construction of the reduced basis $\mathbf{V}\in\mathbb{R}^{N_h\times N}$. In Table~\ref{tab:lifex_RBdim} we report three possible values for the POD tolerance $\varepsilon_{POD}$ and the corresponding RB dimension $N$.

\begin{table}[b!]
	\vspace{-2mm}
	\centering
	\begin{tabular}{|l||c|c|c|} 
		\hline
		POD tolerance $\varepsilon_{POD}$ & $10^{-3}$ & $5\cdot10^{-4}$ & $10^{-4}$ \\
		\hline
		RB dimension $N$ & 39 & 52 & 99 \\
		\hline
	\end{tabular}
	\caption{Cardiac cycle, physiological scenarios. POD tolerances and associated RB dimension for the physiological scenario, when $N_h=18501$.} 
	\label{tab:lifex_RBdim}
	\vspace{-2mm}
\end{table}

Since the input of the encoder function of the DNN-architecture is reshaped into a square matrix (see Remark~\ref{rmk:zero-padded}) and we do not want to introduce too many additional terms when zero-padding, we choose $N$ such that $\sqrt{N+1}\in\mathbb{N}$ (note that the input to $\bm\upsilon_N$ has dimension $N+1$), in this case $N=63$, and build the RB basis by means of the randomized singular value decomposition (SVD). The latter is an efficient, non-deterministic, version of SVD which exploits random sampling to construct a low-dimensional subspace to captures most of the energy of the data matrix, and then manipulates the associated reduced matrix with classical deterministic algorithms, to obtain the desired low-rank approximation.

Once the ROM is built, we perform $n_s'=50$ simulations to collect the reduced data necessary for training the DNNs, namely
\begin{equation}\label{eq:snapshotsDNNs}
\begin{aligned}
&\mathbf{S}_{\bm\rho} = \left[\mathbf{V}^T\mathbf{R}(\mathbf{Vu}_N^{n,(k)}(\bm\mu_\ell),t^{n};\bm\mu_\ell)\right]_{\ell=1,\dots,n_s',n=1,\dots,N_t,k\geq0},\\
&\mathbf{S}_{\bm\iota} = \left[\mathbf{V}^T\partial_\mathbf{u}\mathbf{R}(\mathbf{Vu}_N^{n,(k)}(\bm\mu_\ell),t^{n};\bm\mu_\ell)\mathbf{V}\right]_{\ell=1,\dots,n_s',n=1,\dots,N_t,k\geq0},\\
&\mathbf{S}_{\bm\pi} = \left[\mathbf{V}^T\partial_p\mathbf{R}(\mathbf{Vu}_N^{n,(k)}(\bm\mu_\ell),t^{n};\bm\mu_\ell)\right]_{\ell=1,\dots,n_s',n=1,\dots,N_t,k\geq0},\\
&\mathbf{S}_{\bm\upsilon} = \left[\partial_\mathbf{u}\text{R}^{vol}(\mathbf{Vu}_N^{n,(k)}(\bm\mu_\ell),t^{n};\bm\mu_\ell)\mathbf{V}~|~\text{R}^{vol}(\mathbf{Vu}_N^{n,(k)}(\bm\mu_\ell),t^{n};\bm\mu_\ell)\right ]_{\ell=1,\dots,n_s',n=1,\dots,N_t,k\geq0},
\end{aligned}
\end{equation}
where $\mathbf{S}_{\bm\rho}\in\mathbb{R}^{N\times 1\times N_{train}}$, $\mathbf{S}_{\bm\iota}\in\mathbb{R}^{N\times N\times N_{train}}$, $\mathbf{S}_{\bm\pi}\in\mathbb{R}^{N\times 1\times N_{train}'}$ and $\mathbf{S}_{\bm\upsilon}\in\mathbb{R}^{1\times (N+1)\times N_{train}'}$. Here, $N_{train}$ and $N_{train}'$ denote the total number of snapshots, being $N_{train}>N_{train}'$, since the snapshots for $\mathbf{S}_{\bm\pi}$ and $\mathbf{S}_{\bm\upsilon}$ are collected only during the isovolumetric phases.

A finer computational grid obtained by refining the previous mesh has been also considered to build a second FOM; in this second case, $36704$ elements and $42225$ vertices are used, so that the average cell diameter is equal to 0.0016~m (corresponding to the mesh size commonly used to accurately capture the myocardial displacement with expensive, high-fidelity models \cite{augustin2016anatomically, regazzoni2020cardiac2}). The resulting FOM is character\-ized by $N_h=126675$ degrees of freedom and allows us to assess the performances of Deep-HyROMnet in more realistic situations. In this case, we consider $n_s=15$ samples to build the RB basis $\mathbf{V}$ by means of randomized SVD, choosing $N=63$ as for the coarser mesh. In Table~\ref{tab:lifex_RBdim_fine} we report the RB dimension $N$ obtained for different values of $\varepsilon_{POD}$, showing that the ROM dimension does not increase as $N_h$ becomes larger. Furthermore, we perform $n_s'=50$ ROM simulations in order to build the snapshots matrices (\ref{eq:snapshotsDNNs}) necessary for training the DNNs.

\begin{table}[b!]
	\centering
	\begin{tabular}{|l||c|c|c|} 
		\hline
		POD tolerance $\varepsilon_{POD}$ & $10^{-3}$ & $5\cdot10^{-4}$ & $10^{-4}$ \\
		\hline
		RB dimension $N$ & 37 & 48 & 91 \\
		\hline
	\end{tabular}
	\caption{Cardiac cycle, physiological scenarios. POD tolerances and associated RB dimension for the physiological scenario, when $N_h=126675$.}
	\label{tab:lifex_RBdim_fine}
	\vspace{-2mm}
\end{table}

\begin{remark}
	In order to reduce the computational time required for the training of the DNNs in the case of the finer mesh, we rely on a suitable pre-training strategy \cite{Goodfellow-et-al-2016}, that is, the optimal weights and biases found for the DNNs when $N_h=18501$ are used to initialize the corresponding networks for the larger FOM-dimension.
\end{remark}

Table~\ref{tab:lifex_DEIM_DeepHyROMnet} summarizes the average results obtained on a testing set of 10 input parameters using Deep-HyROMnet for both meshes. Moreover, we report the performances of POD-Galerkin-DEIM built by employing the POD method on the ROM residual snapshots with tolerance $\varepsilon_{DEIM} = 10^{-5}$, corresponding to a DEIM residual basis of dimension $m=1545$. No further speed-ups can be achieved by decreasing the size of the reduced mesh due to convergence issues of the reduced Newton system for some instances of the considered parameters. For $N_h=18501$, Deep-HyROMnet computes a reduced solution in only $16$~s, that is, almost 100 times faster than the reference high-fidelity simulation which requires 27 minutes, whilst yielding an absolute error $\epsilon_{abs}$ on the displacement field of order $O(10^{-2})$. On the other hand, the POD-Galerkin-DEIM ROM, despite being slightly more accurate than Deep-HyROMnet, still requires high computational resources, employing $20$~min to simulate a single heartbeat. When using a finer mesh, Deep-HyROMnet only takes $90$~s to compute the displacement dynamics for a complete heartbeat, against almost $4$~h required by the FOM for the same accuracy level.

\begin{table}[t!]
	\centering
	\begin{tabular}{|l||c|c||c|} 
		\hline
		& ~~~~DEIM-1545~~~~ & Deep-HyROMnet & Deep-HyROMnet \\
		\hline
		$N_h$ & \multicolumn{2}{c||}{$18501$} & $126675$ \\
		\hline
		FOM time & \multicolumn{2}{c||}{27~min} & 3~h 50~min \\
		\hline
		Speed-up & $\times$1.4 & $\times$100 & $\times$150 \\
		\hline
		Avg. CPU time & 20~min & 16~s & 1~min 30~s \\
		\hline
		mean$_{\bm\mu}$ $\epsilon_{abs}(\bm\mu)$ & $3\cdot10^{-3}$ & $3\cdot10^{-2}$ & $7\cdot10^{-2}$ \\
		\hline
		mean$_{\bm\mu}$ $\epsilon_{rel}(\bm\mu)$ & $7\cdot10^{-3}$ & $7\cdot10^{-2}$ & $7\cdot10^{-2}$\\
		\hline
	\end{tabular}
	\caption{Cardiac cycle, physiological scenarios. Computational data related to DEIM-based and DNN-based hyper-ROMs, for $N=63$.}
	\label{tab:lifex_DEIM_DeepHyROMnet}
	\vspace{-2mm}
\end{table}

\begin{figure}[b!]
	\centering
	\includegraphics[width=0.975\textwidth]{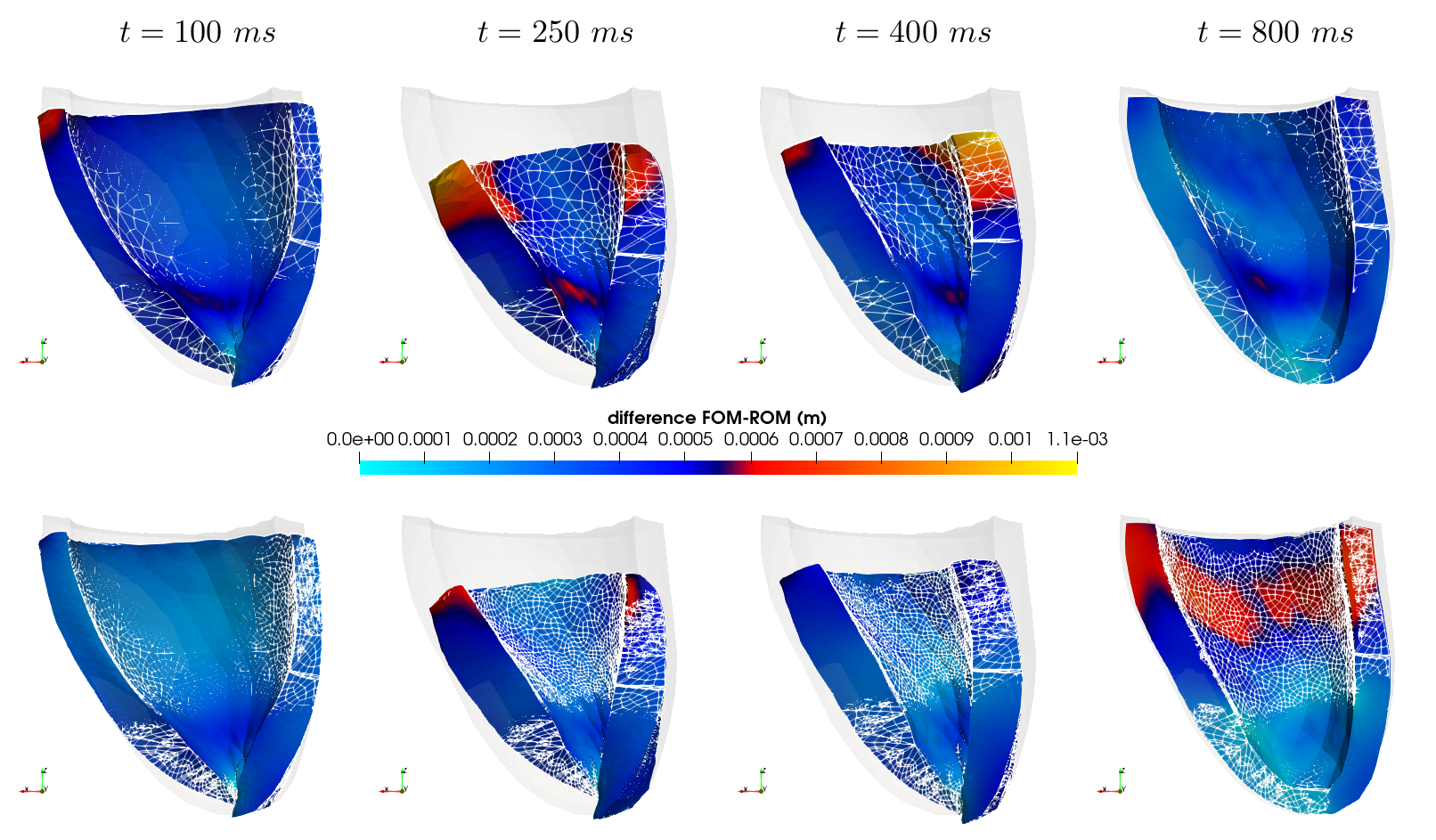}
	\caption{Cardiac cycle, physiological scenarios. FOM (wireframe) and Deep-HyROMnet (colored) displacements for the coarser (top) and the finer (bottom) mesh for $\bm\mu = [34500~\text{Pa},3.00\cdot10^7~\text{Pa}\cdot\text{s}\cdot\text{m}^{-3},55950~\text{Pa}]$.}
	\label{fig:lifex_N63_hyperROM_mu1}
\end{figure}

\begin{figure}[t!]
	\centering
	\includegraphics[width=0.975\textwidth]{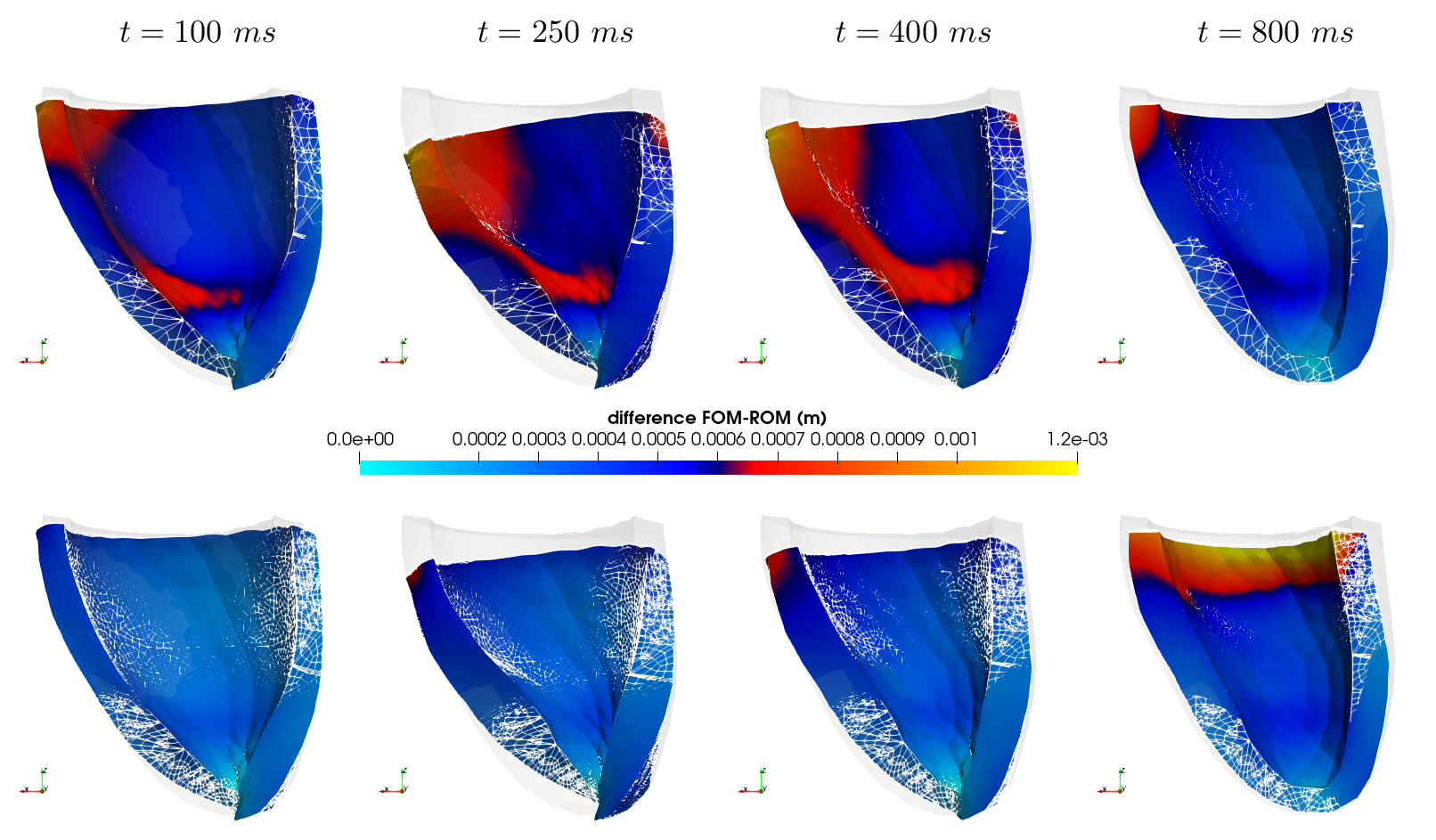}
	\caption{Cardiac cycle, physiological scenarios. FOM (wireframe) and Deep-HyROMnet (colored) displacements for the coarser (top) and the finer (bottom) mesh for $\bm\mu = [58500~\text{Pa},4.16\cdot10^7~\text{Pa}\cdot\text{s}\cdot\text{m}^{-3},49050~\text{Pa}]$.}
	\label{fig:lifex_N63_hyperROM_mu3}
\end{figure}

\begin{figure}[t!]
	\centering
	\includegraphics[width=0.975\textwidth]{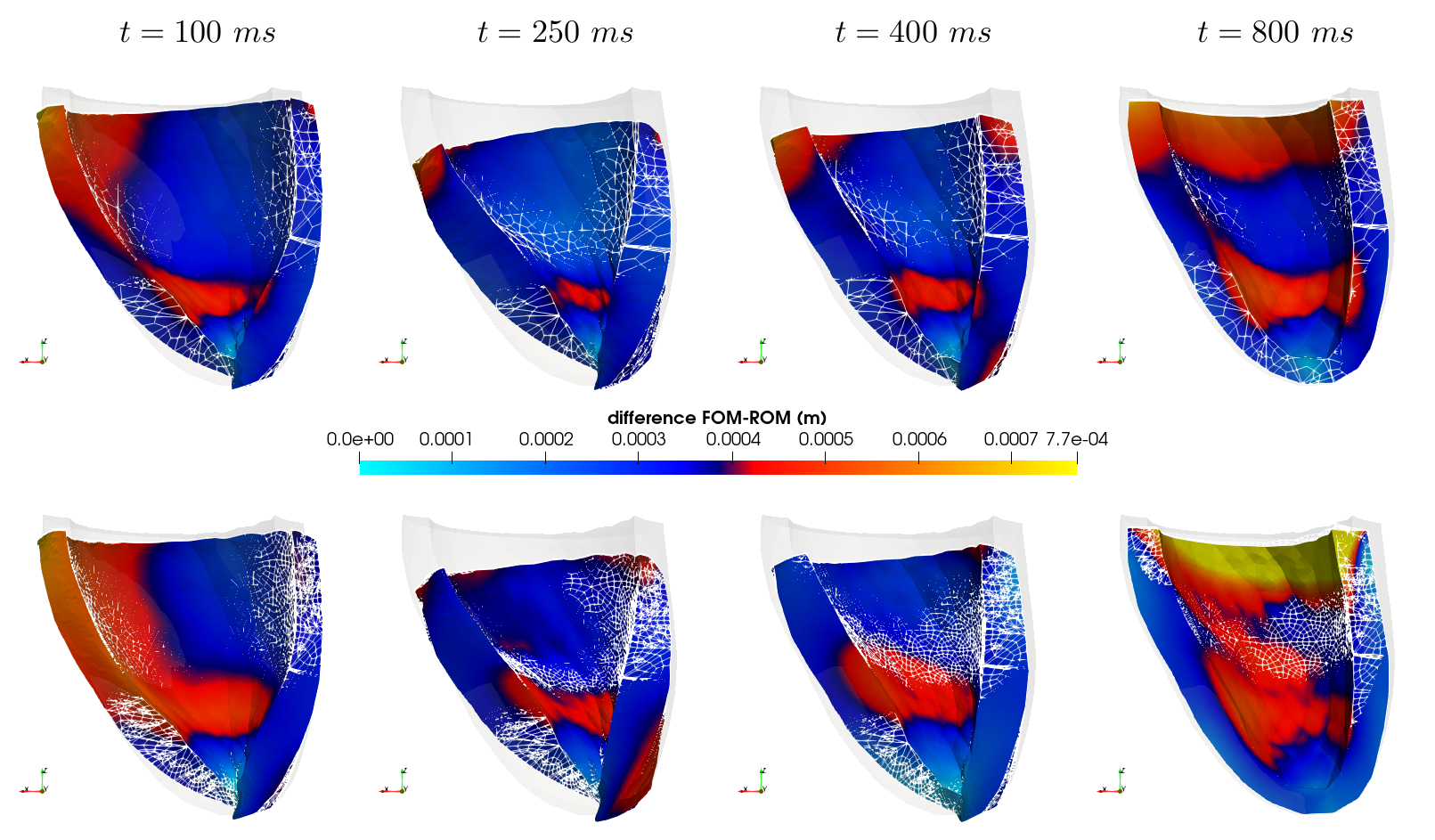}
	\caption{Cardiac cycle, physiological scenarios. FOM (wireframe) and Deep-HyROMnet (colored) displacements for the coarser (top) and the finer (bottom) mesh for $\bm\mu = [66500~\text{Pa},4.20\cdot10^7~\text{Pa}\cdot\text{s}\cdot\text{m}^{-3},57750~\text{Pa}]$.}
	\label{fig:lifex_N63_hyperROM_mu6}
\end{figure}

\begin{figure}[t!]
\vspace{-2mm}
	\centering
		\includegraphics[width=0.72\textwidth]{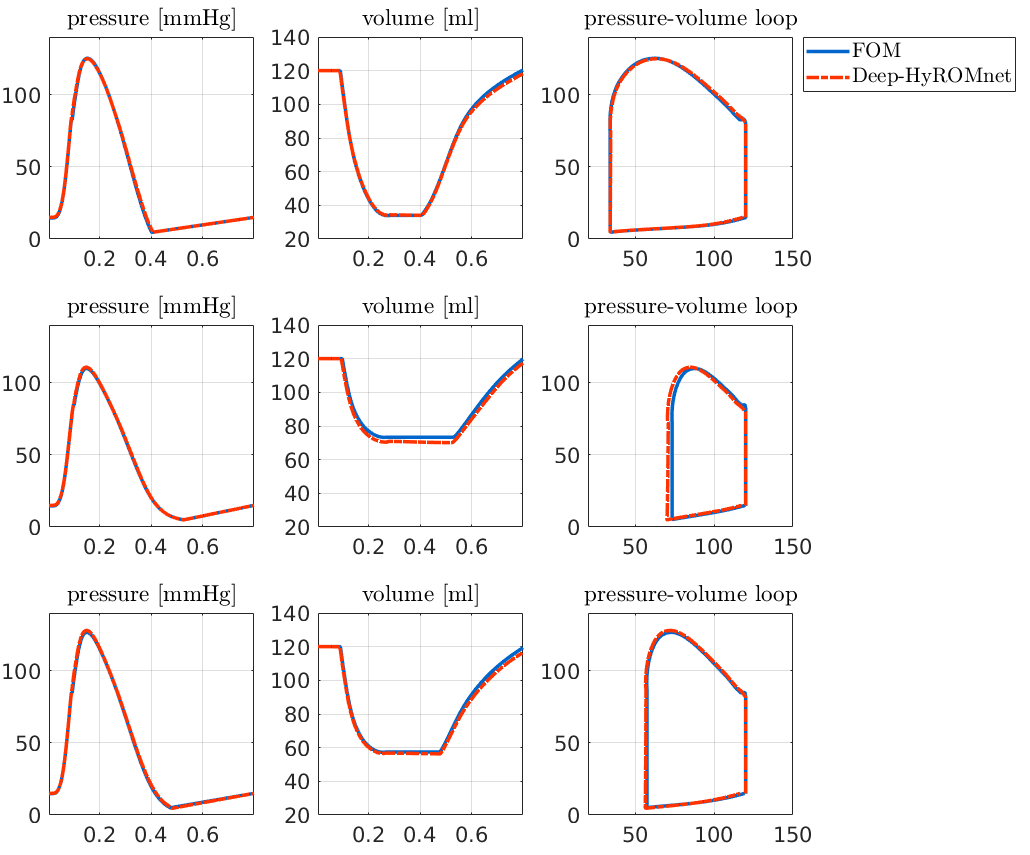}
		\vspace{-2mm}
	\caption{Cardiac cycle, physiological scenarios. Pressures, volumes and pressure-volume relationships, for $\bm\mu = [34500~\text{Pa},3.00\cdot10^7~\text{Pa}\cdot\text{s}\cdot\text{m}^{-3},55950~\text{Pa}]$ (top), $\bm\mu = [58500~\text{Pa},4.16\cdot10^7~\text{Pa}\cdot\text{s}\cdot\text{m}^{-3},49050~\text{Pa}]$ (middle), $\bm\mu = [66500~\text{Pa},4.20\cdot10^7~\text{Pa}\cdot\text{s}\cdot\text{m}^{-3},57750~\text{Pa}]$ (bottom), for $N_h=18501$.}
	\label{fig:pvloop_coarse}
	\vspace{-3mm}
\end{figure}

\begin{figure}[b!]
\vspace{-2mm}
	\centering
	\includegraphics[width=0.72\textwidth]{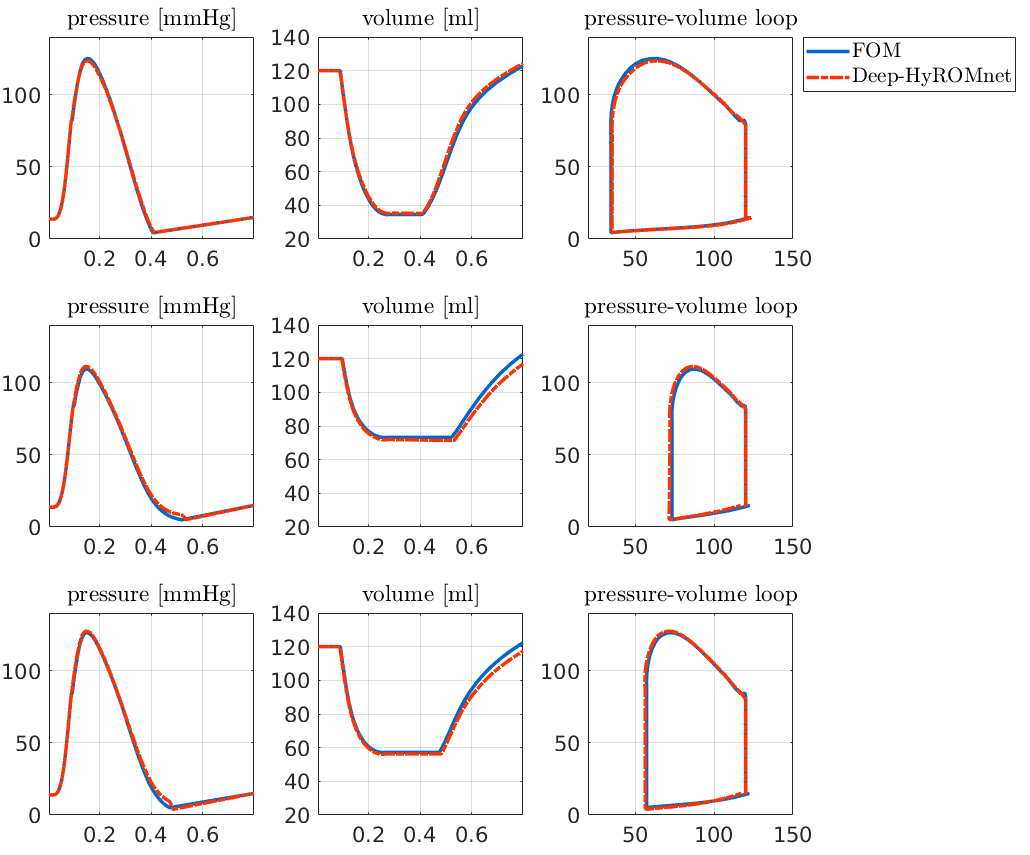}
	\vspace{-2mm}
	\caption{Cardiac cycle, physiological scenarios. Pressures, volumes and pressure-volume relationships, for $\bm\mu = [34500~\text{Pa},3.00\cdot10^7~\text{Pa}\cdot\text{s}\cdot\text{m}^{-3},55950~\text{Pa}]$ (top), $\bm\mu = [58500~\text{Pa},4.16\cdot10^7~\text{Pa}\cdot\text{s}\cdot\text{m}^{-3},49050~\text{Pa}]$ (middle), $\bm\mu = [66500~\text{Pa},4.20\cdot10^7~\text{Pa}\cdot\text{s}\cdot\text{m}^{-3},57750~\text{Pa}]$ (bottom), for $N_h=126675$.}
	\label{fig:pvloop_fine}
\end{figure}

Figures~\ref{fig:lifex_N63_hyperROM_mu1}, \ref{fig:lifex_N63_hyperROM_mu3} and \ref{fig:lifex_N63_hyperROM_mu6} show the Deep-HyROMnet solution for the coarser and the finer mesh computed at different phases of the cardiac cycle, for three different values of the parameter vector. We observe that the pointwise error between the FOM and the Deep-HyROMnet solutions does not increase in time. The corresponding left ventricular pressures and volumes obtained using the FOM and the proposed Deep-HyROM\-net strategy for three input parameters values are reported in Figures~\ref{fig:pvloop_coarse} and \ref{fig:pvloop_fine} for the coarser and the finer test cases, respectively, showing perfect agreement of the reduced outputs of interest with the high-fidelity ones, uniformly on the set of parameter inputs.

\subsection{Pathological scenario} \label{sec:patho}

We now address the solution to the 3D-0D coupled problem in the eventuality that a portion of the cardiac tissue has been affected by myocardial ischemia, that is a reduction of blood supply to the myocardium that may lead to the death of cells in the affected area \cite{griffin2008manual}. In this case, a reduced excitability of the cells and altered ionic currents are observed, as well as inhibited contractility of the tissue \cite{gerbi2018numerical, salvador2021electromechanical}. Reduced order models for cardiac electrophysiology in the case of ischemic necrosis have been formerly considered in, e.g.,  \cite{pagani2018numerical, fresca2020deep}. However, ROMs have never been applied to characterize the mechanical behavior of the myocardium in these scenarios, for varying conditions of the ischemic tissue. Therefore, this is the first time that such a phenomenon is investigated systematically in a broad variety of conditions. 

In this numerical test case, we consider as varying input parameters 
\begin{itemize}
	\item the resistance of the windkessel model $R_p\in[2.5\cdot10^{7},4.5\cdot10^{7}]~\text{Pa}\cdot\text{s}\cdot\text{m}^{-3}$,
	\item the active tension parameter $\widetilde{T}_a\in[4.5\cdot10^{4},6\cdot10^{4}]$~Pa,
	\item and the radius of the ischemic region $r\in[0.2\cdot10^{-3}, 20\cdot10^{-3}]$~m,
\end{itemize}
being them among the most influential parameters associated with the circulation model, the active component of the structural model and the necrotic region, respectively. For the computational domain we employ the hexahedral mesh with $6167$ vertices reported in Figure~\ref{fig:LV_loaded_BCs} (center), so that the FOM obtained using $\mathbb{Q}_1$-FE has dimension $N_h=18501$.

\begin{figure}[b!]
	\centering
	\vspace{-0.1cm}
	\includegraphics[width=0.8\textwidth]{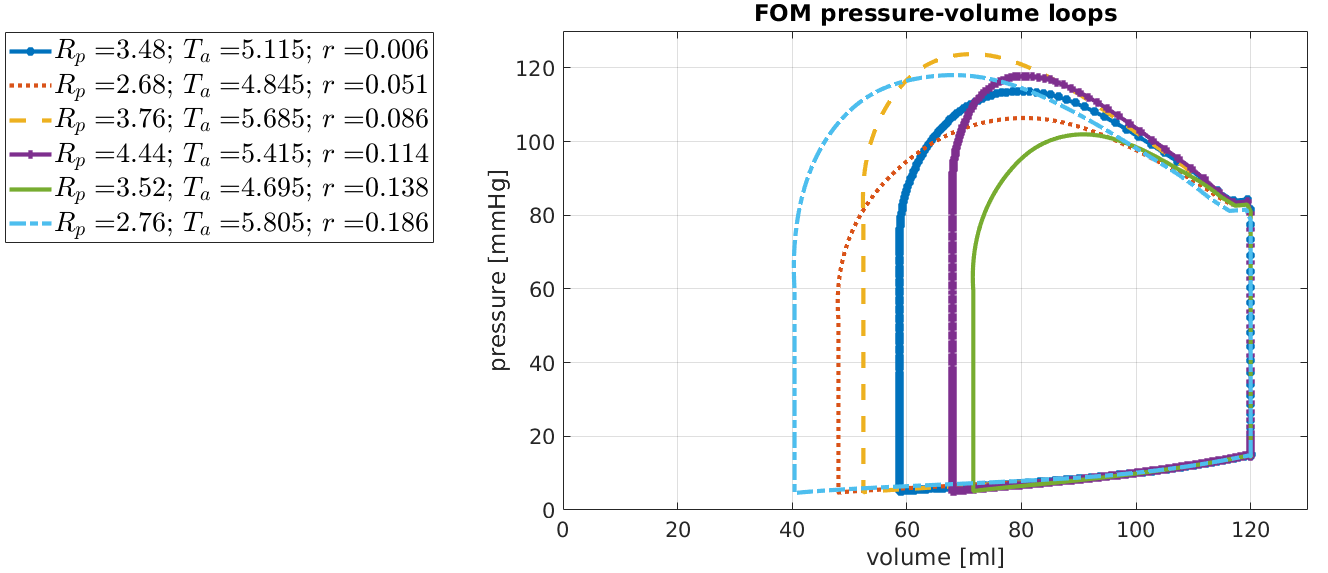}
	\vspace{-0.2cm}
	\caption{Cardiac cycle, pathological scenarios. Pressure-volume loops computed using the FOM with different parameter instances.}
	\label{fig:lifex_ischemic_FOM_pv}
\end{figure}

Figure~\ref{fig:lifex_ischemic_FOM_pv} shows the pressure-volume loops obtained for six different values of the input vector $\bm\mu = [R_p,T_a,r]\in\mathcal{P}$, highlighting how parameter variations may have a great impact on outputs of interest. In particular, we observe that the end systolic volume ranges from $40$~ml to $72$~ml, so that the ejection fraction ($EF$), that is the volumetric portion of blood ejected from the ventricle with each contraction, reduces from $67\%$ to $40\%$. On the other hand, the maximum values of the blood pressure goes from $102$~mmHg to almost $124$~mmHg, thus influencing the slope of the end systolic pressure-volume relationship which provides an index of myocardial contractility \cite{sato1998espvr}. All these indicators are useful in clinical practice \cite{doyle2013left, burkhoff2005assessment, shoucri2015end}. However, in order to gain more knowledge about the relationships between model inputs and outputs of interest, sensitivity analysis studies have to be performed.

The reduced basis $\mathbf{V}$ is built by collecting high-fidelity solution snapshots for $n_s=20$ parameter samples and performing randomized SVD using $N=80$. In fact, a higher dimension of the RB basis with respect to the physiological scenario is required, possibly due to the presence of the ischemia, as highlighted in Table~\ref{tab:lifex_ischemic_RBdim}. Indeed, the presence of the scar region ultimately makes the parameters-to-solution map more involved, affecting the behavior of the solution in a more pronounced way and thus requiring a higher dimension of the basis if a global linear subspace has to be used to approximate the whole solution manifold with sufficient accuracy.

\begin{table}[t!]
	\vspace{2mm}
	\centering
	\begin{tabular}{|l||c|c|c|} 
		\hline
		POD tolerance $\varepsilon_{POD}$ & $10^{-3}$ & $5\cdot10^{-4}$ & $10^{-4}$ \\
		\hline
		RB dimension $N$ & 72 & 100 & 184 \\
		\hline
	\end{tabular}
	\caption{Cardiac cycle, pathological scenarios. POD tolerances and associated RB dimension for the pathological scenario, when $N_h=18501$.}
	\label{tab:lifex_ischemic_RBdim}
	\vspace{2mm}
\end{table}

Then, $n_s'=50$ POD-Galerkin ROM simulations are performed to collect the reduced nonlinear data (\ref{eq:snapshotsDNNs}) and the DNNs are trained. Table~\ref{tab:lifex_ischemic_DEIM_DeepHyROMnet} reports the results obtained using Deep-HyROMnet, where the average is computed over 20 testing parameters.

\begin{table}[t!]
	\vspace{-2mm}
	\centering
	\begin{tabular}{|l||c|c|} 
		\hline
		& DEIM-3000 & Deep-HyROMnet \\
		\hline
		Speed-up & $\times$1.2 & $\times$94 \\
		\hline
		Avg. CPU time & 22~min & 17~s \\
		\hline
		mean$_{\bm\mu}$ $\epsilon_{abs}(\bm\mu)$ & $6\cdot10^{-3}$ & $5\cdot10^{-2}$ \\
		\hline
		mean$_{\bm\mu}$ $\epsilon_{rel}(\bm\mu)$ & $2\cdot10^{-3}$ & $1\cdot10^{-1}$\\
		\hline
	\end{tabular}
	\caption{Cardiac cycle, pathological scenarios. Computational data related to DEIM-based and DNN-based hyper-ROMs, for $N=80$.}
	\label{tab:lifex_ischemic_DEIM_DeepHyROMnet}
	\vspace{-0.25cm}
\end{table}

\begin{figure}[b!]
\vspace{-0.25cm}
	\centering
	\includegraphics[width=0.75\textwidth]{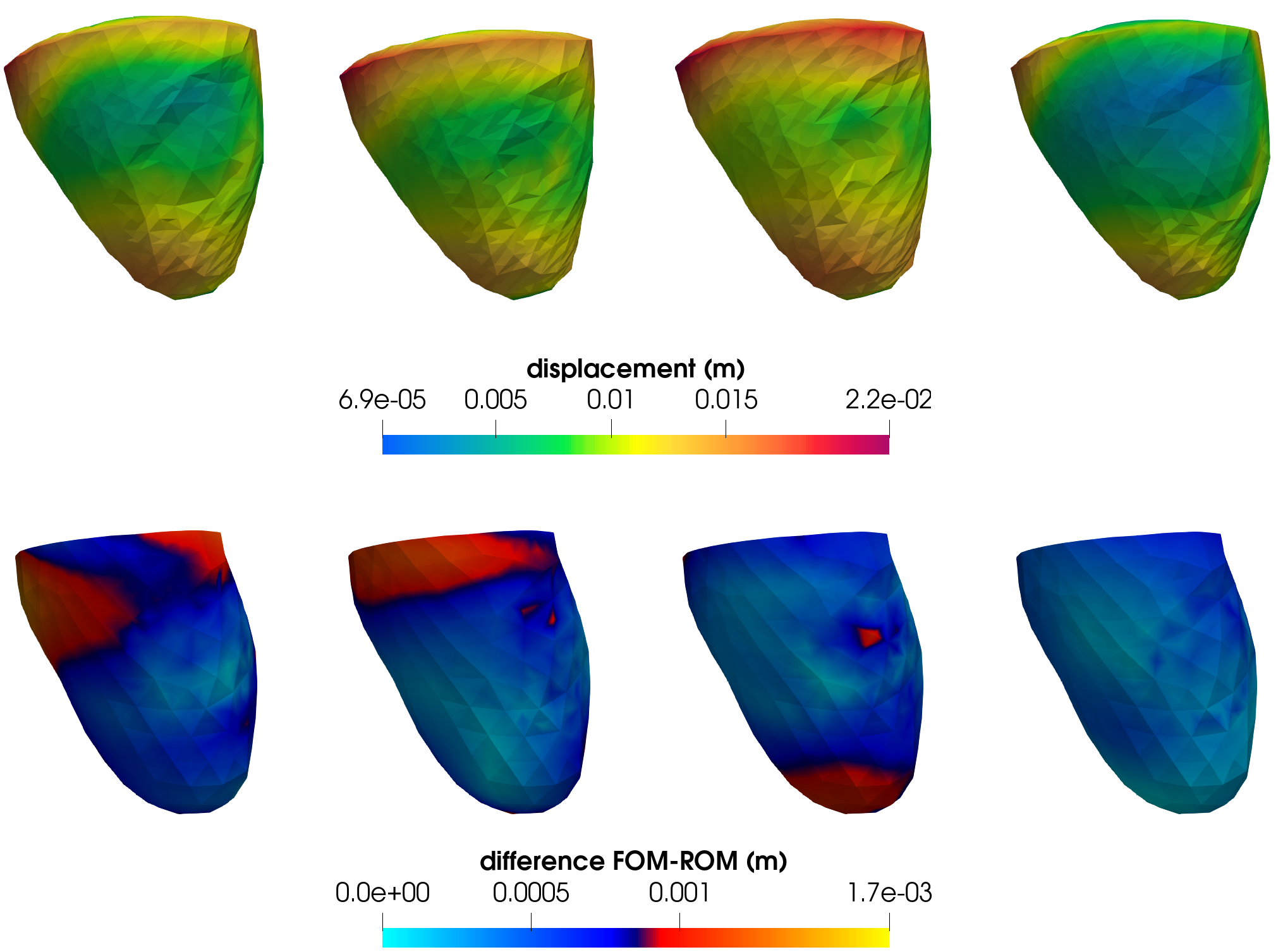}
	\caption{Cardiac cycle, pathological scenarios. Deep-HyROMnet deformation (top) and error (bottom) at time $t=0.25$~s for different values of the parameter $\bm\mu$ (from left to right).}
	\label{fig:lifex_ischemic_DeepHyROMnet}
	\vspace{-0.25cm}
\end{figure}

\begin{figure}[t!]
	\centering
	\includegraphics[width=0.85\textwidth]{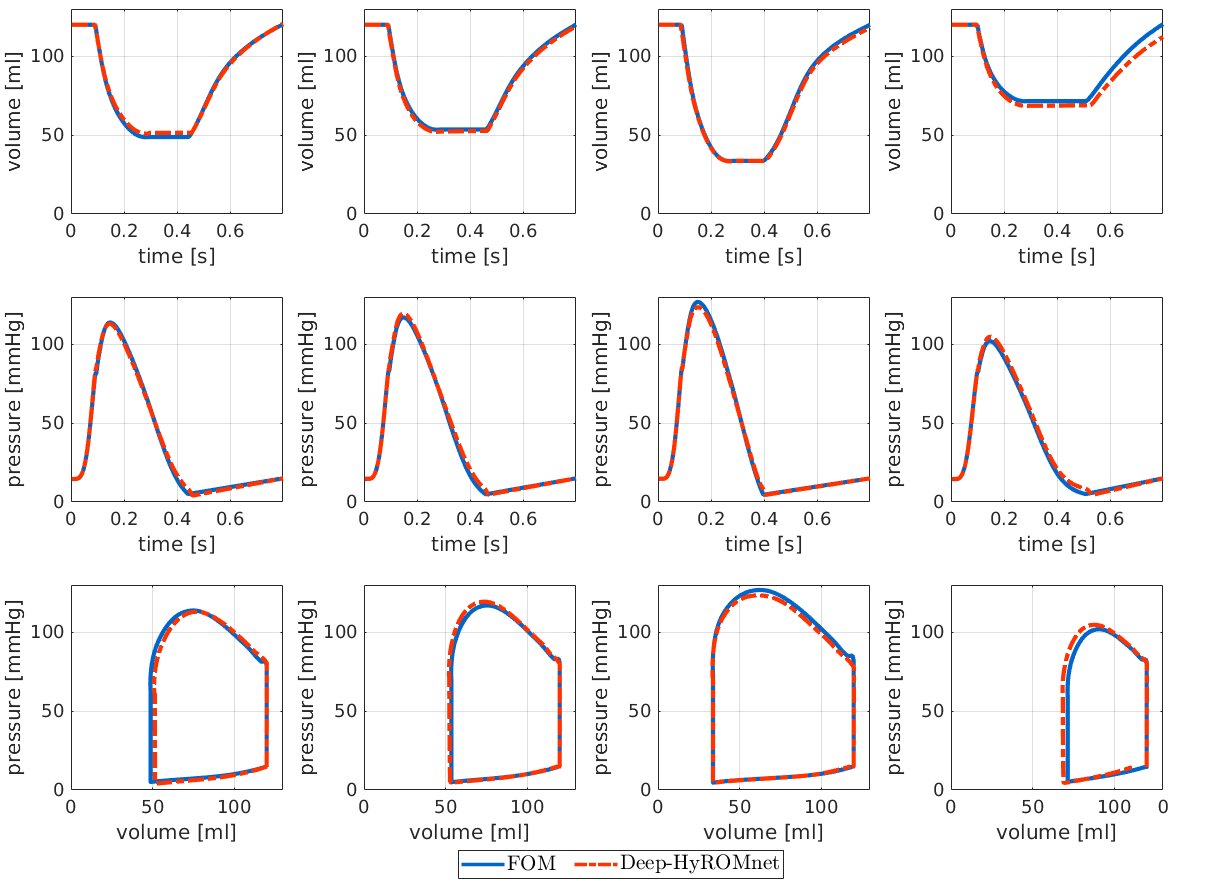}
	\caption{Cardiac cycle, pathological scenarios. Pressure, volumes and pressure-volume relationships for different values of the parameter $\bm\mu$ (from left to right).}
	\label{fig:lifex_ischemic_DeepHyROMnet_pvloops}
\end{figure}

As in the physiological scenario, Deep-HyROMnet requires less than $17$~s to compute a whole heartbeat, so that it is almost $100$ times faster than the FOM, which for the same task to be achieved requires almost 26 minutes. It is worth mentioning that relying on POD-Galerkin-DEIM ROMs led to negligible speed-ups (e.g., only $1.2$ times faster that the FOM), thus making the development of the Deep-HyROMnet ROM necessary to efficiently address the solution to the problem under investigation. For what concerns the accuracy of the hyper-ROM, the absolute displacement error $\epsilon_{abs}$ is around $5\cdot10^{-2}$. Although more accurate results can be obtained with classical hyper-reduction techniques, a good approximation of the outputs of interest is obtained using our DNN-based ROM. In particular, the error between the FOM and Deep-HyROMnet on the $EF$, computed over the testing set, is less than $3\%$. Figure~\ref{fig:lifex_ischemic_DeepHyROMnet} reports few examples of hyper-ROM displacement and pointwise error with respect to the FOM at time $t=0.25$~s, when the ventricle in fully contracted, while the corresponding pressure-volume loops are shown in Figure~\ref{fig:lifex_ischemic_DeepHyROMnet_pvloops}.

\section{Application to Forward Uncertainty Quantification}\label{sec:forwardUQ}

To conclude, in this section we address the repeated evaluation of the inputs-to-solution map in both physiological and pathological scenarios by means of the Deep-HyROMnet ROMs developed in Sections~\ref{sec:physio} and \ref{sec:patho}, in order to gain some useful knowledge about the impact of the model parameters on selected output quantities. With this aim, let us consider as output quantities of interest:
\begin{itemize}
	\item the ejection fraction, that represents the amount of blood pumped at each heartbeat, and is defined as
	\begin{equation*}
	EF = \frac{EDV-ESV}{EDV},
	\end{equation*}
	where $EDV$ and $ESV$ denote the end-diastolic and the end-systolic volumes, respectively;
	\item the maximal rate of change in pressure 
	\begin{align*}
	dP/dt_{max} = \underset{t\in(0,T)}{\max}\left(\frac{d{p}_{LV}(t;\bm\mu)}{dt}\right)\approx \underset{n=1,\dots,N_t}{\max}\left( \frac{p_{LV}^n(\bm\mu)-p_{LV}^{n-1}(\bm\mu)}{\Delta t} \right),
	\end{align*}
	which is a common indicator of cardiac contractility. 
\end{itemize}
These choices are motivated by the fact that both $EF$ and $dP/dt_{max}$ are commonly used mechanical biomarkers. Nonetheless, since Deep-HyROMnet computes the whole displacement at each time instance, any additional output, such as, e.g., the wall thickening, the end-systolic pressure or the longitudinal fractional shortening \cite{levrero2020sensitivity, campos2020uncertainty}, can be considered online without the need to rebuild the ROM. This is a distinguishing feature of the proposed reduction technique, compared to recent frameworks addressing NN-based approximation of quantities of interest, without taking into account the approximation of the field variables involved in the output evaluations \cite{regazzoni2021realtime}.

\newpage

For what concerns the varying parameters, we consider:
\begin{itemize}
	\item the resistance of the windkessel model $R_p\in[2.5\cdot10^{7},4.5\cdot10^{7}]~\text{Pa}\cdot\text{s}\cdot\text{m}^{-3}$,
	\item the active tension parameter $\widetilde{T}_a\in[4.5\cdot10^{4},6\cdot10^{4}]$~Pa, and
	\item the radius of the ischemic region $r\in[0.2\cdot10^{-3}, 20\cdot10^{-3}]$~m,
\end{itemize}
so that we always assume $K=5\cdot10^4$~Pa as online value for the bulk modulus in the healthy case as well.

The following results are obtained by performing $500$ hyper-ROM simulations in physiological scenar\-i\-os, i.e. for $\bm\mu = [R_p,\widetilde{T}_a]$, and $1000$ in pathological ones, i.e. for $\bm\mu = [R_p,\widetilde{T}_a, r]$, taking into account an hexahedral computational mesh of a patient-specific left ventricle (see Figure~\ref{fig:LV_loaded_BCs}) with $6167$ vertices. We recall that in this case the underlying FOM dimension is $N_h=18501$ and that less than $17$~s are required by Deep-HyROMnet to compute the problem solution for each new parameter instance, thus entailing less than $7$ hours of CPU time on a PC desktop computer with 3.70GHz Intel Core i5-9600K CPU and 16GB RAM. Performing these studies using the FOM would have required $27$ days of computations, which becomes almost $240$ days if the finer computational grid with $42225$ vertices has to be considered (reducing to only $38$ hours when employing Deep-HyROMnet).

\begin{figure}[b!]
	\centering
	\includegraphics[width=0.9\textwidth]{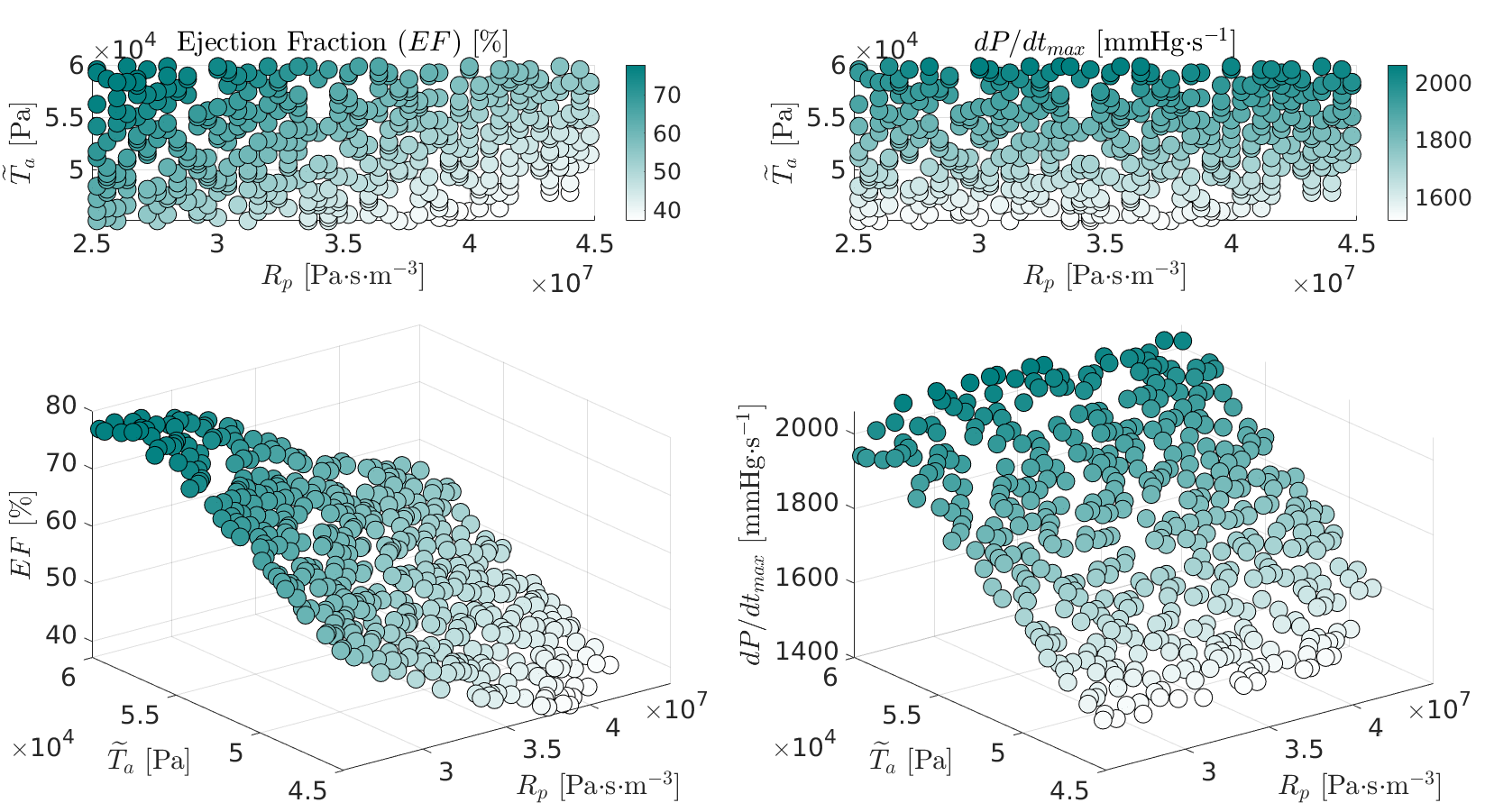}
	\caption{Cardiac cycle, physiological scenarios. Scatter plots of the $EF$ (left) and the $dP/dt_{max}$ (right) for $500$ different parameters. The colormap represents the values of the output of interest.}
	\label{fig:3D_physio}
\end{figure}

Concerning the outcomes of the healthy scenarios, we observe that both the resistance $R_p$ of the two-element windkessel model and the active stress parameter $\widetilde{T}_a$ have a great impact on the $EF$, as shown in Figure~\ref{fig:3D_physio} (left). Both variables are, in fact, associated with the systolic phase of the cardiac cycle: larger values of the maximum active tension lead to a greater contraction of the myocardial tissue, whereas higher values of $R_p$ correspond to a lower amount of blood that the ventricle is able to pump during ejection (phase~\ref{it:ejection} of the heartbeat, see Figure~\ref{fig:Wiggers}). As a consequence, they both affect the $ESV$ without substantially changing the $EDV$. In particular, given a fixed value of $R_p$, the $EF$ increases as $\widetilde{T}_a$ becomes higher; on the other way round, when $\widetilde{T}_a$ is fixed, the $EF$ decreases as the resistance of the circulation model is increased. As an example, the minimum value $EF=37\%$ corresponds to $\bm\mu=[4\cdot10^7~\text{Pa}\cdot\text{s}\cdot\text{m}^{-3},4.605\cdot10^4~\text{Pa}]$, that is when $R_p$ and $\widetilde{T}_a$ are closed to their upper and lower bounds, respectively; its maximum value $EF=78\%$ is obtained instead for $\bm\mu=[2.6\cdot10^7~\text{Pa}\cdot\text{s}\cdot\text{m}^{-3},5.835\cdot10^4~\text{Pa}]$. On the other hand, from Figure~\ref{fig:3D_physio} (right), we can conclude that the maximal rate of change in pressure is proportional to the active stress, going from $1763~\text{mmHg}\cdot\text{s}^{-1}$ to $1998~\text{mmHg}\cdot\text{s}^{-1}$ as $\widetilde{T}_a$ is increased from $4.515\cdot10^4$~Pa to its maximum value $5.985\cdot10^4$~Pa, whilst we observe that $R_p$ has almost no influence on $dP/dt_{max}$.

\begin{figure}[b!]
\vspace{-2mm}
	\centering
	\begin{subfigure}{0.49\textwidth}
		\centering
		\includegraphics[width=\textwidth]{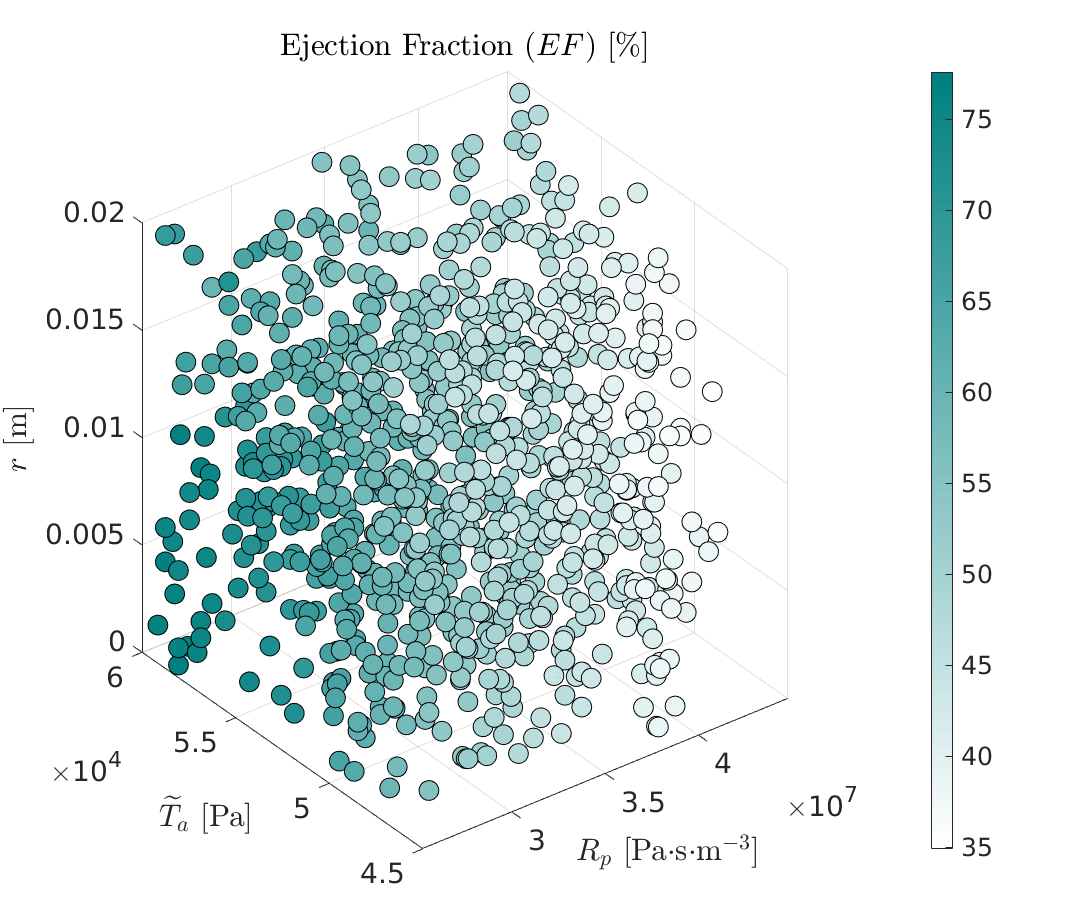}
		\caption{$EF$}
	\end{subfigure}
	\begin{subfigure}{0.49\textwidth}
		\centering
		\includegraphics[width=\textwidth]{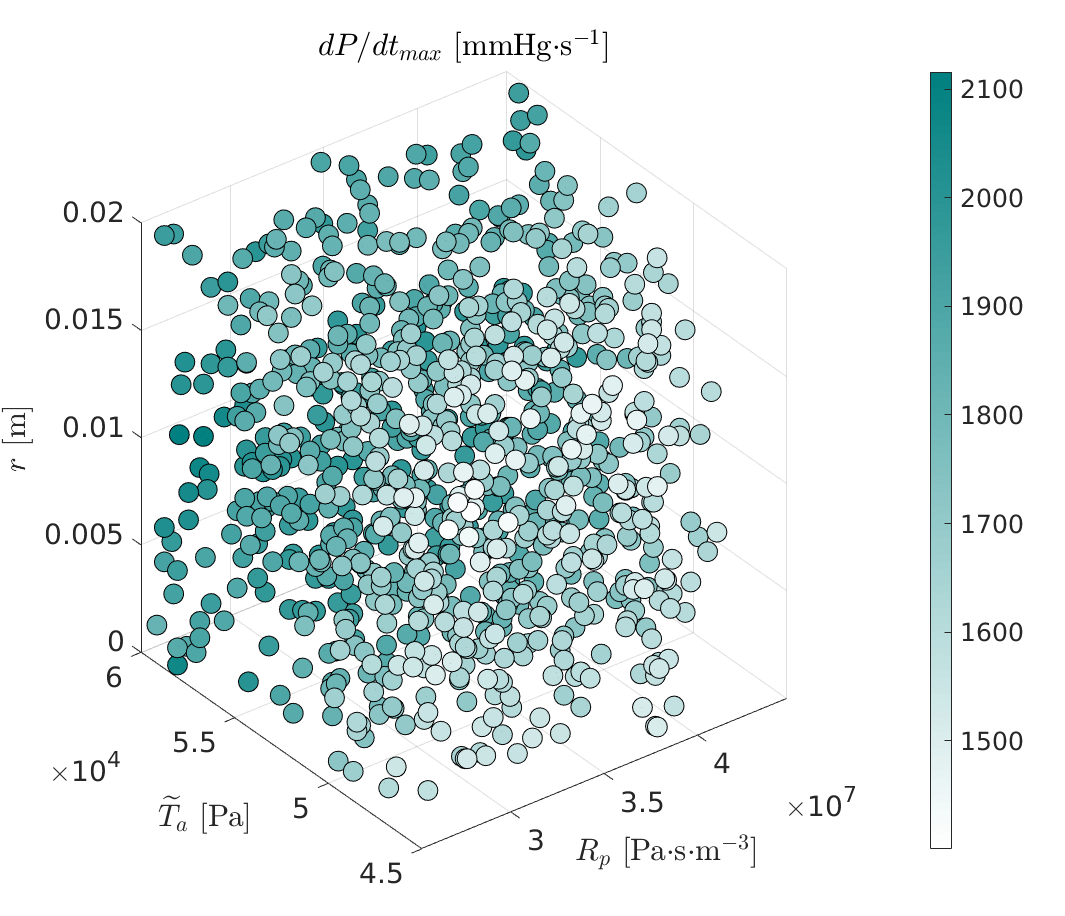}
		\caption{$dP/dt_{max}$}
	\end{subfigure}
	\vspace{-2mm}
	\caption{Cardiac cycle, pathological scenarios. Scatter plots of the $EF$ (left) and the $dP/dt_{max}$ (right) in the $xyz$-plane ($x=R_p$, $y=\widetilde{T}_a$ and $z=r$) for $1000$ different parameters. The colormap represents the values of the outputs of interest.} \vspace{2mm}
	\label{fig:3D_patho}
	\vspace{-2mm}
\end{figure}

Assessing the way input variations affect the considered outputs of interest in the pathological scenarios becomes more involved due to the presence of an additional parameter, and to the fact that no activation of the cardiac myocytes is assumed inside the necrotic region $\mathcal{B}(\mathbf{X}_c,r)$, being $\mathbf{X}_c\in\Omega_0$ a fixed point inside the myocardium. In Figure~\ref{fig:3D_patho} we report the scatter plots of $EF$ and $dP/dt_{max}$, where in the $x$-,$y$- and $z$-axis are reported $R_p$, $\widetilde{T}_a$ and the radius $r$, respectively, while the colors of the data points encode the value of the outputs $EF$ and $dP/dt_{max}$.

\begin{figure}[b!]
\vspace{-2mm}
	\centering
	\includegraphics[width=\textwidth]{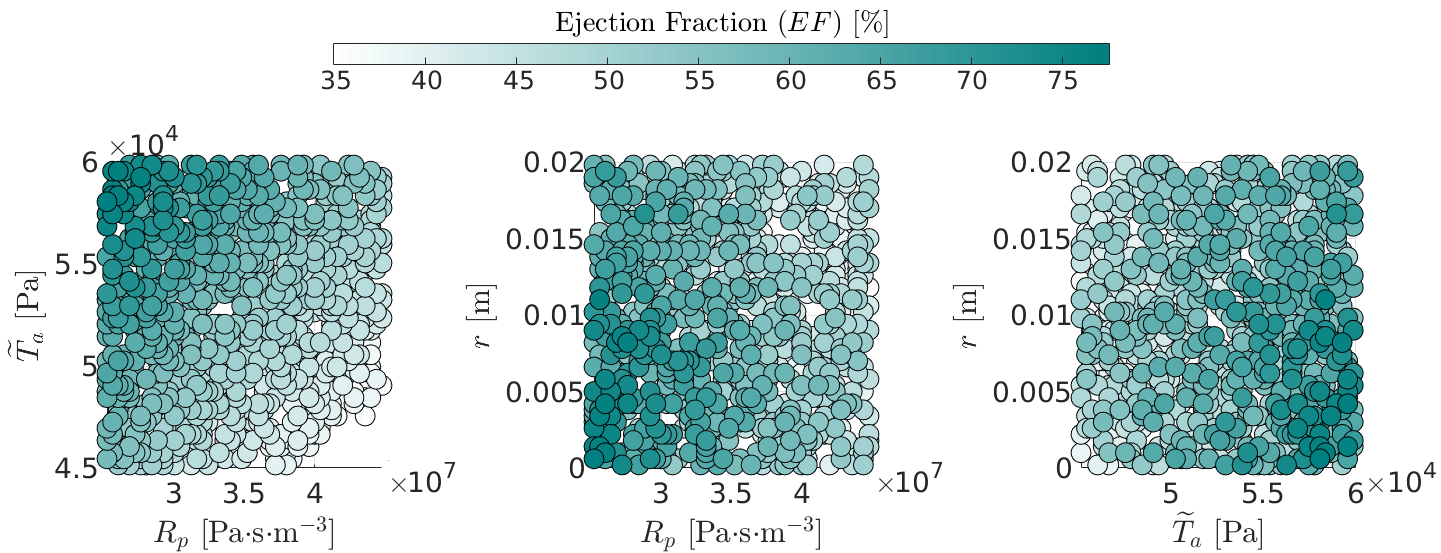}
	\vspace{-2mm}
	\caption{Cardiac cycle, pathological scenarios. 2D-views of the scatter plots of the $EF$ for $1000$ different parameters. The colormap represents the values of the output of interest.}
	\label{fig:2D_patho_EF}
	\vspace{-2mm}
\end{figure}

Regarding the interaction between the maximum active tension $\widetilde{T}_a$ and the windkessel resistance $R_p$ on their influence on the $EF$ (Figure~\ref{fig:2D_patho_EF}, left), we can draw similar conclusions to the healthy case. To give few examples, the lower values of $EF=34,9\%$ and $35.6\%$ are computed for the parameters $\bm\mu=[4.48\cdot10^7~\text{Pa}\cdot\text{s}\cdot\text{m}^{-3},4.905\cdot10^4~\text{Pa},0.0118\cdot\text{m}^{-3}]$ and $\bm\mu=[4.36\cdot10^7~\text{Pa}\cdot\text{s}\cdot\text{m}^{-3},4.755\cdot10^4~\text{Pa},0.0066\cdot\text{m}^{-3}]$, respectively, whereas the higher values $EF=76.8\%$ and $77.6\%$ are associated with the inputs $\bm\mu=[2.60\cdot10^7~\text{Pa}\cdot\text{s}\cdot\text{m}^{-3},5.835\cdot10^4~\text{Pa},0.0018\cdot\text{m}^{-3}]$ and $\bm\mu=[2.56\cdot10^7~\text{Pa}\cdot\text{s}\cdot\text{m}^{-3},5.835\cdot10^4~\text{Pa},0.011\cdot\text{m}^{-3}]$. On the other hand, the influence of $r$ on the $EF$ is more difficult to ascertain from the analysis of the scatter plots.

\begin{figure}[t!]
\vspace{-2mm}
	\centering
	\includegraphics[width=\textwidth]{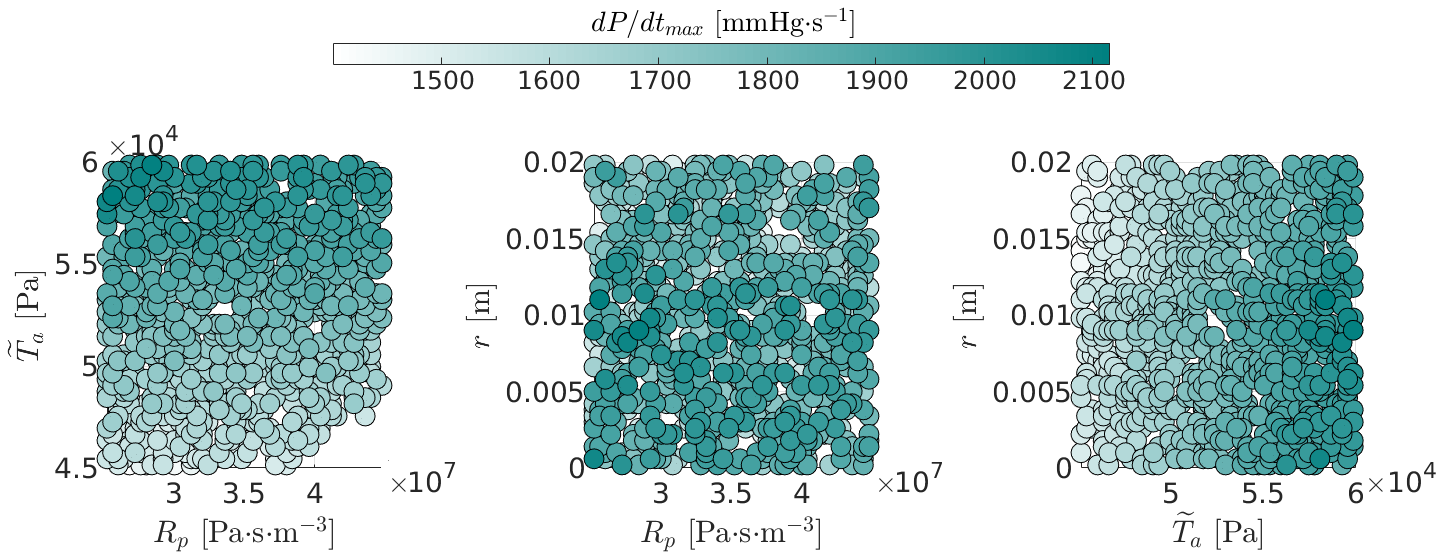}
	\caption{Cardiac cycle, pathological scenarios. 2D-views of the scatter plots of the $dP/dt_{max}$ for $1000$ different parameters. The colormap represents the values of the output of interest.}
	\label{fig:2D_patho_dPdt}
	\vspace{-2mm}
\end{figure}

Finally, from the 2D-views of the scatter plots reported in Figure~\ref{fig:2D_patho_dPdt}, we can assume that variations of both $R_p$ and $r$ have almost to effect on the maximal rate of change of pressure $dP/dt_{max}$, and that $\widetilde{T}_a$ is the most influential parameter between those considered.

To conclude, we have observed that the maximum value of the active tension $\widetilde{T}_a$ has great influence on $EF$ and $dP/dt_{max}$, both in the physiological and in the pathological tests considered. This fact is compliant with the results of sensitivity analysis conducted in \cite{campos2020uncertainty} for the healthy left ventricle in a quasi-static assumption. Furthermore, the resistance $R_p$ of the circulation model, associated with the ejection phase of ventricular systole, influences the values the $EF$, as well as the size of the necrosis.


\section{Conclusion}\label{sec:conclusion}

In this work, we have applied our new physics-based, (hyper-)reduced order modeling strategy, named \textit{Deep-HyROMnet}, for the accurate and efficient approximation of nonlinear elastodynamics problems arising in cardiac mechanics. This method combines POD for the construction of a reduced basis, Galerkin projection over the low-dimensional subspace spanned by these basis functions and DNNs for efficiently handle the nonlinear reduced operators. We proved that Deep-HyROMnet is able to obtain an extremely good approximation of the displacement field, as well as pressure and volume useful to compute key scalar cardiac outputs, with few reduced dofs, while achieving considerable speed-ups thanks to the approximation of the nonlinear terms by means of a DNN architecture. Our method is suitable for a range of scenarios in which classical projection-based ROMs  would require high computational costs.

In particular, we have shown how Deep-HyROMnets outperform POD-Galerkin-DEIM ROMs in terms of computational speed-up and allow to address the efficient solution to cardiac mechanics problems coupled with a lumped-parameter model for blood circulation, both in physiological and pathological scenarios. Preliminary results of forward uncertainty quantification carried out on a patient-specific left ventricle allowed to gain some useful knowledge about the impact of the model parameters on possible output quantities of interest. In this context, we have observed that the active tension has great influence on both the ejection fraction and the maximal rate of change in pressure; other parameters, such as the resistance of the circulation model and the size of the necrosis, showed instead higher influence on the ejection fraction only.

By providing a reliable and computationally efficient reduction procedure, our model can be successfully used to address the solution of multi-query problems, such as, e.g., forward uncertainty quantification and parameter estimation. However, further investigations are needed to assess the performance of the proposed reduction strategy on even more involved scenarios. Remarkable examples include {\em (i)} the use of even finer computational meshes and smaller time steps; {\em (ii)} the introduction of a surrogate model for the computation of space- and time-dependent active tension, thus taking into account the activation of cardiac myocytes at different time instants; {\em (iii)}   more and/or different  input parameters according to sensitivity analysis carried on the simulation of the whole cardiac cycle.

\section*{Acknowledgements}

The authors have been supported by the ERC Advanced Grant iHEART, “An integrated heart model for the simulation of the cardiac function”, 2017-2022, P.I. A. Quarteroni (ERC2016AdG, project ID: 740132). 
The authors gratefully acknowledge Dr. R. Piersanti, Dr. F. Regazzoni and Dr. M. Salvador (MOX, Politecnico di Milano) for their useful remarks and discussions regarding the coupled electromechanical model and the calibration of the active force generation surrogate model, as well as Dr. P. Africa (MOX, Politecnico di Milano) for his kind support while using some methods implemented in the in-house Finite Element library \texttt{life$^\texttt{x}$}.

\bibliography{Biblio_Cicci_et_al_3.bib}
\bibliographystyle{ieeetr}


\begin{appendices}
	
	
	\section{The POD technique}\label{appendix:POD}
	Given $n_s<N_h$ randomly sampled instances of the parameter $\bm\mu\in\mathcal{P}$, define the snapshots matrix 
	\begin{equation*}
	\mathbf{S} = \left[ \mathbf{u}_h(t^1;\bm\mu_1)~|~\dots~|~\mathbf{u}_h(t^{N_t};\bm\mu_1)~|~\dots~|~\mathbf{u}_h(t^1;\bm\mu_{n_s})~|~\dots~|~\mathbf{u}_h(t^{N_t};\bm\mu_{n_s}) \right]\in\mathbb{R}^{N_h\times n_s},
	\end{equation*}
	where $\mathbf{u}_h(t^n;\bm\mu_\ell)$ are FOM solutions computed for different values of $\bm\mu_\ell$. The RB basis $\mathbf{V}\in\mathbb{R}^{N_h\times N}$ is obtained by performing the singular valued decomposition of $\mathbf{S}$
	\begin{equation*}
	\mathbf{S}=\mathbf{U\Sigma Z}^T,
	\end{equation*}
	and collecting the first $N$ columns of $\mathbf{U}\in\mathbb{R}^{N_h\times n_s}$, corresponding to the first $N$ left singular vectors. This yields an orthonormal basis that, among all $N$-dimensional orthonormal basis, minimizes the least square error of the snapshot reconstruction, that it
	\begin{align*}
	\lVert \mathbf{S} - \mathbf{V}\mathbf{V}^T\mathbf{S} \rVert_F^2 = \underset{\{\mathbf{W}\in\mathbb{R}^{N_h\times N}~|~\mathbf{W}^T\mathbf{W}=\mathbf{I}\}}{\min} \lVert \mathbf{S} - \mathbf{W}\mathbf{W}^T\mathbf{S} \rVert_F^2 = \sum_{i=N+1}^r \sigma_i^2,
	\end{align*}
	where $\lVert \cdot \rVert_F$ is the Frobenius norm and $\sigma_1\geq\dots\geq\sigma_{r}\geq0$ are the singular values of $\mathbf{S}$, being rank$(\mathbf{S})=r\leq N_h\land n_s$. Thus, the singular values provide a quantitative criteria for choosing the size $N$, which is typically computed as the minimum integer satisfying the condition
	\begin{equation}\label{eq:N-POD}
	\frac{\sum_{i=1}^{N}\sigma_i^2}{\sum_{i=1}^{r}\sigma_i^2} \geq 1-\varepsilon_{POD}^2
	\end{equation}
	for a given tolerance $\varepsilon_{POD}>0$. The POD technique is summarized in Algorithm \ref{alg:POD}.
	
	\begin{algorithm}
		\caption{Proper orthogonal decomposition (POD)}\label{alg:POD}
		INPUT: snapshots matrix $\mathbf{S}\in\mathbb{R}^{N_h\times n_s}$, target tolerance $\varepsilon_{POD}>0$\\
		OUTPUT: RB basis $\mathbf{V}\in\mathbb{R}^{N_h\times N}$, with $N = N(\varepsilon_{POD})$
		\begin{algorithmic}[1]
			\STATE Perform singular value decomposition of $\mathbf{S}$: $\mathbf{S=U\Sigma Z}^T$
			\STATE Select the basis dimension $N$ as the minimum integer fulfilling condition (\ref{eq:N-POD})
			\STATE Construct $\mathbf{V}$ selecting the first $N$ columns of $\mathbf{U}$
		\end{algorithmic}
	\end{algorithm}

	
	\section{DL-ROM-based neural network}\label{appendix:DeepHyROMnet}	
	For the sake of completeness, we briefly describe the DNN-based approximation of the reduced residual vector, that is
	\begin{equation*}
	{\bm\rho}_N(\bm\mu,t^n,k) \approx\mathbf{R}_N(\mathbf{V}\mathbf{u}_N^{n,(k)}(\bm\mu),t^n;\bm\mu)\in\mathbb{R}^N.
	\end{equation*}
	We point out that, by defining the transformation
	\begin{equation*}
	vec\colon\mathbb{R}^{N\times N}\rightarrow\mathbb{R}^{N^2}, \quad vec(\mathbf{J}_N(\mathbf{V}\mathbf{u}_N^{n,(k)}(\bm\mu),t^n;\bm\mu)) = \mathbf{j}_N(\mathbf{V}\mathbf{u}_N^{n,(k)}(\bm\mu),t^n;\bm\mu),
	\end{equation*}
	which consists in stacking the columns of $\mathbf{J}_N(\mathbf{V}\mathbf{u}_N^{n,(k)}(\bm\mu),t^n;\bm\mu)$ in a vector of dimension $N^2$, we can applied the DL-ROM technique described for the residual on the Jacobian matrix as well, thus obtaining 
	\begin{equation*}
	\widetilde{\bm\iota}_N(\bm\mu,t^n,k)\approx \mathbf{j}_N(\mathbf{V}\mathbf{u}_N^{n,(k)}(\bm\mu),t^n;\bm\mu)\in\mathbb{R}^{N^2}.
	\end{equation*}
	Finally, the $vec$ operation is reverted to obtain ${\bm\iota}_N(\bm\mu,t^n,k) = vec^{-1}(\widetilde{\bm\iota}_N(\bm\mu,t^n,k))$.
	
\begin{figure}
		\centering
		\includegraphics[width=\textwidth]{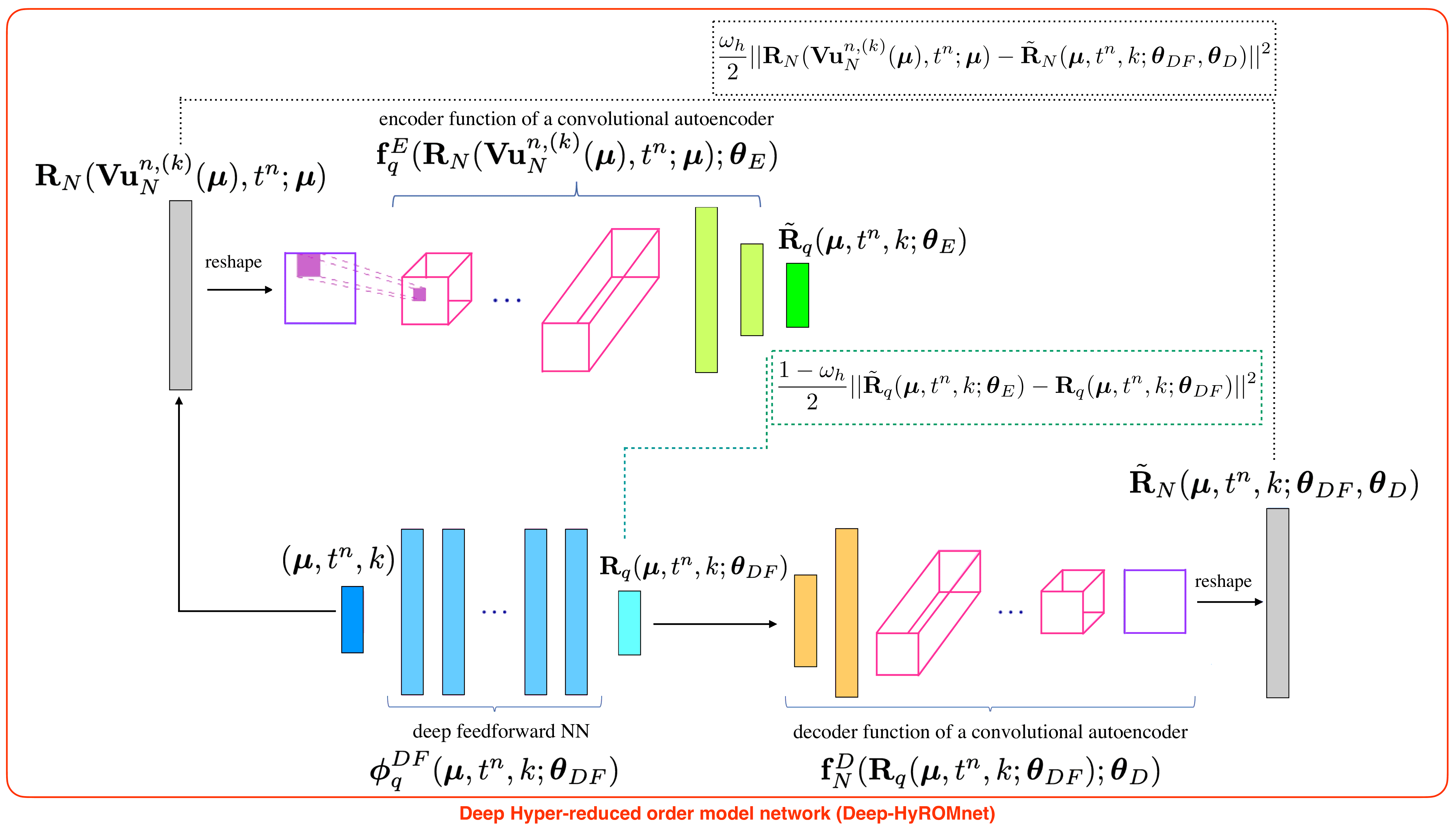}
		\caption{DNN architecture used during the training phase for the reduced residual vector.}
		\label{fig:DNN}
	\end{figure}

	The DL-ROM approximation of the ROM residual $\mathbf{R}_N(\mathbf{V}\mathbf{u}_N^{n,(k)}(\bm\mu),t^n;\bm\mu)$ takes the form
	\begin{equation*}
	\bm\rho_N(\bm\mu,t^n,k) = \widetilde{\mathbf{R}}_N(\bm\mu,t^n,k;\bm\theta_{DF},\bm\theta_{D}) = \mathbf{f}^D_N(\bm\phi_q^{DF}(\bm\mu,t^n,k;\bm\theta_{DF});\bm\theta_{D})
	\end{equation*}
	where
	\begin{itemize}
		\item $\bm\phi_q^{DF}(\cdot~;\bm\theta_{DF})\colon\mathbb{R}^{P+2}\rightarrow\mathbb{R}^q$ such that 
		\begin{equation*}
		\bm\phi_q^{DF}(\bm\mu,t^n,k;\bm\theta_{DF}) = \mathbf{R}_q(\bm\mu,t^n,k;\bm\theta_{DF})
		\end{equation*} is a deep feedforward neural network (DFNN), where $\bm\theta_{DF}$ denotes the vector of parameters, collecting all the corresponding weights and biases of each layer, and $q$ is as close as possible to the input size $P+2$;
		\item $\mathbf{f}^D_N(\cdot~;\bm\theta_{D})\colon\mathbb{R}^q\rightarrow\mathbb{R}^N$ such that
		\begin{equation*}
		\mathbf{f}^D_N(\mathbf{R}_q(\bm\mu,t^n,k;\bm\theta_{DF});\bm\theta_{D}) = \widetilde{\mathbf{R}}_N(\bm\mu,t^n,k;\bm\theta_{DF},\bm\theta_{D})
		\end{equation*}
		is the decoder function of a convolutional autoencoder (CAE), depending upon the vector $\bm\theta_{D}$ of weights and biases.
	\end{itemize}
	The encoder function of the CAE is exploited, during the training stage only, to map the reduced residual $\mathbf{R}_N(\mathbf{V}\mathbf{u}_N^{n,(k)}(\bm\mu),t^n;\bm\mu)$ associated to $(\bm\mu,t^n,k)$ onto a low-dimensional representation
	\begin{equation*}
	\mathbf{f}^E_q(\mathbf{R}_N(\mathbf{V}\mathbf{u}_N^{n,(k)}(\bm\mu),t^n;\bm\mu);\bm\theta_{E}) = \widetilde{\mathbf{R}}_q(\bm\mu,t^n,k;\bm\theta_{E}),
	\end{equation*}
	where $\mathbf{f}^E_q(\cdot~;\bm\theta_{E})\colon\mathbb{R}^N\rightarrow\mathbb{R}^q$ denotes the encoder function and $\bm\theta_{E}$ is the corresponding vector of parameters. The architecture used during training is reported in Figure~\ref{fig:DNN}, whereas, during the testing phase, the encoder function $\mathbf{f}^E_q$ is discarded.
	
	\begin{remark}\label{rmk:zero-padded}
		The input of the encoder function. i.e. $\mathbf{R}_N$, is reshaped into a square matrix by rewriting its elements in row-major order, thus obtaining $\mathbf{R}_N^{reshape}\in\mathbb{R}^{\sqrt{N}\times \sqrt{N}}$. If $N$ is not a square, the input is zero-padded as explained in \cite{Goodfellow-et-al-2016}, and the additional elements are subsequently discarded.	
	\end{remark}
	
	Let
	\begin{equation*}	
	\mathbf{S}_{\bm\rho} = \left[\mathbf{R}_N(\mathbf{V}\mathbf{u}_N^{n,(k)}(\bm\mu_\ell),t^n;\bm\mu_\ell)\right]_{\ell=1,\dots,n_s',n=1,\dots,N_t,k\geq0}\in\mathbb{R}^{N\times N_{train}},
	\end{equation*}
	with $N_{train}=n_s'N_tN_k$, be the reduced residual snapshots matrix collecting column-wise ROM residuals computed for $n_s'$ sampled parameters $\bm\mu_\ell\in\mathcal{P}$, at different time instances $t^1,\dots,t^{N_t}$ and for each Newton iteration $k\geq0$. Moreover, we define the parameter matrix of the corresponding triples as
	\begin{equation*}
	\mathbf{M} = \left[\left(\bm\mu_\ell,t^n,k\right)\right]_{\ell=1,\dots,n_s',n=1,\dots,N_t,k\geq0}\in\mathbb{R}^{(P+2)\times N_{train}}.
	\end{equation*}
	The training stage consists in solving the following optimization problem in the weights variable $\bm\theta = (\bm\theta_{E},\bm\theta_{DF},\bm\theta_{D})$:
	\begin{equation*}
	\underset{\bm\theta}{\min}~\mathcal{J}(\bm\theta) = \underset{\bm\theta}{\min } \dfrac{1}{N_{train}}\sum_{\ell=1}^{n_s'}\sum_{n=1}^{N_t}\sum_{k=0}^{N_k}\mathcal{L}(\bm\mu_\ell,t^n,k;\bm\theta),
	\end{equation*}
	where
	\begin{equation}\label{eq:loss}
	\begin{aligned}
	\mathcal{L}(\bm\mu_\ell,t^n,k;\bm\theta) = &\dfrac{\omega_h}{2} \lVert \mathbf{R}_N(\mathbf{V}\mathbf{u}_N^{n,(k)}(\bm\mu_\ell),t^n;\bm\mu_\ell) -  \widetilde{\mathbf{R}}_N(\bm\mu_\ell,t^n,k;\bm\theta_{DF},\bm\theta_{D}) \rVert^2 \\
	& + \dfrac{1-\omega_h}{2} \lVert \widetilde{\mathbf{R}}_q(\bm\mu_\ell,t^n,k;\bm\theta_{E}) - \mathbf{R}_q(\bm\mu_\ell,t^n,k;\bm\theta_{DF}) \rVert^2,
	\end{aligned}
	\end{equation}
	with $\omega_h\in[0,1]$. For further details on the training and testing stages, as well as the corresponding algorithms, we refer to \cite{cicci2021DeepHyROMnet}.
	
	\section{Reference values for the 3D-0D coupled simulations}
	Here we report the reference values used throughout this work for the mechanics and circulation models (if not otherwise specified).
	\begin{table}[h]
		\vspace{-2mm}
		\centering
		\begin{tabular}{|llcc|} 
			\hline
			\textbf{Name} & \textbf{Parameter} & \textbf{Value} & \textbf{Unit} \\
			\hline
			\multicolumn{4}{|c|}{Cardiac mechanics} \\
			\hline
			Tissue density & $\rho_0$ & $10^3$ & $\text{kg}\cdot\text{m}^{-3}$ \\	
			Robin boundary condition & $K_{\perp}$ & $2\cdot10^5$ & $\text{Pa}\cdot \text{m}^{-1}$ \\
			Robin boundary condition & $K_{\parallel}$ & $2\cdot10^4$ & $\text{Pa}\cdot \text{m}^{-1}$ \\
			Robin boundary condition & $C_{\perp}$ & $2\cdot10^4$ & $\text{Pa}\cdot \text{s}\cdot \text{m}^{-1}$ \\
			Robin boundary condition & $C_{\parallel}$ & $2\cdot10^3$ & $\text{Pa}\cdot \text{s}\cdot \text{m}^{-1}$ \\
			\hline
			\multicolumn{4}{|c|}{Passive myocardial tissue} \\
			\hline
			Hyperelastic parameter & $b_f$ & 8 & \\
			Hyperelastic parameter & $b_s$ & 6 & \\
			Hyperelastic parameter & $b_n$ & 3 & \\
			Hyperelastic parameter & $b_{fs}$ & 12 & \\
			Hyperelastic parameter & $b_{fn}, b_{sn}$ & 3 & \\
			Material stiffness & $C$ & 880 & Pa \\
			Bulk modulus & $K$ & $5\cdot10^4$ & Pa \\
			\hline
			\multicolumn{4}{|c|}{Active myocardial tissue} \\
			\hline
			Maximum active tension & $\widetilde{T}_a$ & $5\cdot10^{4}$ & Pa \\
			Fiber angle & $\bm\alpha^{epi}$ & $-60$ & deg \\
			Fiber angle & $\bm\alpha^{endo}$ & $60$ & deg \\
			Fiber angle & $\bm\beta^{epi}$ & $20$ & deg \\
			Fiber angle & $\bm\beta^{endo}$ & $-20$ & deg \\
			\hline
		\end{tabular}
		\caption{Reference values of the input parameters to the 3D mechanics model.}
		\label{tab:reference_values_3D}
		\vspace{-2mm}
	\end{table}

	\begin{table}[h]
		\vspace{-2mm}
		\centering
		\begin{tabular}{|llcc|} 
			\hline
			\textbf{Name} & \textbf{Parameter} & \textbf{Value} & \textbf{Unit} \\
			\hline
			\multicolumn{4}{|c|}{Circulation} \\
			\hline
			Capacitance & $C_p$ & $4.5\cdot10^{9}$ & $\text{m}^{-3}\cdot\text{Pa}^{-1}$\\
			Resistance & $R_p$ & $3.5\cdot10^{7}$ & $\text{Pa}\cdot\text{s}\cdot\text{m}^{-3}$\\
			End-diastolic pressure & $p_{ED}$ & $15$ & mmHg\\
			Aortic valve opening pressure & $p_{AVO}$ & $82,50$ & mmHg\\
			Mitral valve opening pressure & $p_{MVO}$ & $5$ & mmHg \\
			\hline
		\end{tabular}
		\caption{Reference values of the input parameters to the 0D circulation model.}
		\label{tab:reference_values_0D}
		\vspace{2mm}
	\end{table}

\end{appendices}

\end{document}